%% file: main.tex
\documentclass[final,onefignum,onetabnum]{siamonline190516}

\pdfoutput=1

\usepackage[section]{placeins} % For '\input{file}'
\input{preamble}
\begin{document}

\title{An extended range of stable flux reconstruction schemes on quadrilaterals for various polynomial bases\thanks{\corresponding{Will Trojak (\email{w.trojak@imperial.ac.uk})} \funding{None to declare.}}}
\headers{FR on quadrilaterals}{W. Trojak et al.}

\author{Will Trojak\thanks{Department of Aeronautics, Imperial College London, South Kensington, London, SW7 2AZ} \and Rob Watson\thanks{Aeronautical and Automotive Engineering, Loughborough University, Loughborough, LE11 3TU} \and Peter Vincent\footnotemark[2].}

\maketitle

\begin{abstract}
    An extended range of energy stable flux reconstruction schemes, developed using a summation-by-parts approach, is presented on quadrilateral elements for various sets of polynomial bases. For the maximal order bases, a new set of correction functions which result in stable schemes is found. However, for a range of orders it is shown that only a single correction function can be cast as a tensor-product. Subsequently, correction functions are identified using a generalised analytic framework that results in stable schemes for total order and approximate Euclidean order polynomial bases on quadrilaterals --- which have not previously been explored in the context of flux reconstruction. It is shown that the approximate Euclidean order basis can provide similar numerical accuracy as the maximal order basis but with fewer points per element, and thus lower cost.
\end{abstract}


%\linenumbers

\input{introduction}
\input{prelim}
\input{linear}
\input{erfr_quad_mo}
\input{erfr_quad_to}
\input{erfr_quad_eo_near}
\input{numerical}

\section{Conclusions}\label{sec:conclusions}
    Three sets of linearly stable high-order flux reconstruction schemes on quadrilateral elements have been presented. These three sets were formed for the maximal order, total order, and approximate Euclidean order polynomial bases. For the maximal order bases, it has been shown that the previously used tensor product of one-dimensional correction functions do not form part of this set, except for the DG correction functions themselves. Through numerical experimentation with the different bases, it was shown that the Euclidean order basis had similar performance to the maximal order basis, despite using fewer points, and was also significantly more isotropic than the total order basis. 
    This result is consistent with previous observations made when using similar bases for polynomial interpolation. Future work will go on to investigate the utility of Euclidean basis polynomials in FR for real world non-linear problems.

\section*{Acknowledgements}\label{sec:ack}
    WT would like to thank Nick Trefethen for his useful discussions. 

%\section*{Data Availability}
%The data that support the findings of this study are available from the corresponding author upon reasonable request.

\bibliographystyle{plainnat}
\bibliography{reference}


\begin{appendices}
\end{appendices}

% Show the list of todo's in the document.  Needed to avoid stupid warnings/errors when using the todo package
%\todos

\end{document}

%% file: preamble.tex
\usepackage{algorithmic}
\usepackage[caption=false]{subfig}
\ifpdf
  \DeclareGraphicsExtensions{.eps,.pdf,.png,.jpg}
\else
  \DeclareGraphicsExtensions{.eps}
\fi

\usepackage{amssymb}
\usepackage{amsmath}
\usepackage{amsfonts}
\usepackage{upgreek}
\usepackage{siunitx}
\usepackage{listings}
\usepackage{multirow}
\usepackage{multicol}
\usepackage{bigdelim}
\usepackage{color}
\usepackage{stmaryrd}
\usepackage{soul}
\usepackage[symbol]{footmisc}
\usepackage{booktabs}
\usepackage[section]{placeins} % For '\input{file}'
\usepackage[titletoc]{appendix}
\usepackage{calc}
\usepackage[capitalize]{cleveref}
\usepackage[sort&compress,square,numbers]{natbib}
\usepackage{xspace}

\usepackage{graphicx}
\usepackage{graphics}
\usepackage{float}
\usepackage{wrapfig}
\usepackage[percent]{overpic}
\usepackage{varwidth}
\usepackage{epstopdf}
\usepackage{adjustbox}
\usepackage{tikz}
\usepackage{pgfplots}
\usetikzlibrary{arrows,matrix,positioning,fit}
\usetikzlibrary{shapes,positioning}
\usetikzlibrary{backgrounds}
\pgfplotsset{compat=1.14}
\usepgfplotslibrary{colorbrewer}
\usepgfplotslibrary{patchplots}
\usepgfplotslibrary[colorbrewer]
\usetikzlibrary{pgfplots.colorbrewer}
\usetikzlibrary[pgfplots.colorbrewer]
\usepgfplotslibrary{fillbetween}
\usepackage{tikz-layers}
\usepgfplotslibrary{units}
\usetikzlibrary{spy}
\usepackage{pgfplotstable}
\graphicspath{ {figs/} }
\usetikzlibrary{external}

\newcommand\hb[1]{\hat{\mathbf{#1}}}

\newcommand\tb[1]{\tilde{\mathbf{#1}}}

\newcommand\px[2]{\frac{\partial #1}{\partial {#2}}}

\newcommand\dx[2]{\frac{\mathrm{d} #1}{\mathrm{d} #2}}

\newcommand{\half}{\frac{1}{2}}

\newlength\myheight
\newlength\mydepth
\settototalheight\myheight{Xygp}
\newsiamremark{remark}{Remark}
\newsiamremark{hypothesis}{Hypothesis}
\crefname{hypothesis}{Hypothesis}{Hypotheses}
\newsiamthm{claim}{Claim}

\DeclareRobustCommand\onedot{\futurelet\@let@token\@onedot}
\def\@onedot{\ifx\@let@token.\else.\null\fi\xspace}

\def\ie{\emph{i.e}\onedot}

\makeatother

%% file: introduction.tex
\section{Introduction}\label{sec:intro}
    The high-order flux reconstruction (FR) method of \citet{Huynh2007} is an efficient and versatile method for approximating the solution of time dependent partial differential equations. Many works have explored a range of the analytical characteristics of FR in one-dimension~\citep{Vincent2010,Vincent2011,Vincent2015,Ranocha2016,Asthana2017}, but fewer works have studied FR as it is applied to quadrilaterals. Two works which have explored quadrilaterals and the stability of the method when correction functions are formed of a tensor-product of one-dimensional corrections functions are \citet{Sheshadri2015} and \citet{Cicchino2021}. Of these, \citet{Sheshadri2015} was able to construct a stability proof using surface terms which are not reconcilable with the analytical approaches of \citet{Vincent2015} and \citet{Ranocha2016}. In the original study by \citet{Huynh2007} and in a later work by \citet{Trojak2020}, the properties of the FR method were explored using Fourier analysis on quadrilaterals --- and stark differences were observed in the numerical properties of the method when the correction function was changed. Again, both papers made use of a tensor-product of one-dimensional schemes. In the context of implicit large eddy simulation (ILES), \citet{Vermeire2016} has shown that aliasing errors can be greatly affected by the correction function when using a tensor-product of the one-dimensional schemes defined by \citet{Vincent2015}.
    
    The definition of stable FR schemes on quadrilaterals has been entirely limited to these tensor-product schemes, whereas on triangles wide sets of stable FR schemes have been defined --- notably the sets of \citet{Castonguay2011} and of \citet{Williams2013}. More recently, the summation-by-parts (SBP) methods have gained significant research attention due to their utility in the analysis of methods. Using the SBP framework, \citet{Ranocha2016} was able to define an extended set of one-dimensional stable FR schemes. More recently still, \citet{Trojak2022a} have made use of this method to extend the set of stable FR schemes on triangles. 
    
    Within the literature on finite elements, it has been common across many applications for the approximation space on a quadrilateral elements to make use of a maximal order polynomial basis. For example, a first order maximal order basis would include the terms $1$, $x$, $y$, and $xy$. This does fit naturally with the element, but other choices are also compatible. In two works, Trefethen~\citep{Trefethen2017,Trefethen2017b} explored the effect of using other bases when approximating functions, and showed that the so-called Euclidean basis often performs nearly as well as a maximal order basis, but at a lower computational cost. However, this work did overlook one advantageous aspect of the maximal order basis: on quadrilaterals it allows for operators to be decomposed to utilise the tensor product for improved computational efficiency~\citet{Swirydowicz2019,Trojak2022}.
    
    In this work, we will make use of the SBP methods set out by \citet{Trojak2022a} to produce an extended range of stable FR methods on quadrilaterals with a maximal order polynomial basis. This SBP approach will then be generalised to produce analogous sets of stable schemes for alternative bases. With these sets of stable schemes defined, we will go on to investigate the isotropy of the different bases to determine the potential suitability of lower-cost bases. Consequently, this work is structured with the preliminaries given in \cref{sec:prelim} and the key requirements for stability and symmetry defined in \cref{sec:linear}. Then, in \cref{sec:quad_mo}, an extended range of stable FR methods for the maximal order basis is presented and the stability of tensor-product constructions investigated. In \cref{sec:quad_to,sec:quad_neo}, additional sets of stable FR schemes are defined on two alternative polynomial bases, namely the total order basis and an approximate Euclidean order basis. In \cref{sec:numerical}, some numerical tests are presented for the three bases, and finally, in \cref{sec:conclusions}, various conclusions are drawn.

%% file: prelim.tex
\section{Preliminaries}\label{sec:prelim}
\subsection{Flux Reconstruction}
    The flux reconstruction (FR) scheme was first introduced by \citet{Huynh2007} and has been applied to several element topologies and to both advection and advection-diffusion systems~\citep{Huynh2009,Castonguay2012}. To give a brief introduction to the FR method here, we will consider the advection equation in one dimension:
    \begin{equation}\label{eq:advect1d}
        \px{u}{t} + \px{f}{x} = 0,\quad \mathrm{for} \quad u(x,t):K\times \mathbb{R}_+ \mapsto\mathbb{R}, \quad \mathrm{and} \quad f(u):\mathbb{R}\mapsto\mathbb{R}.
    \end{equation}
    The FR algorithm makes use of a sub-division of the domain $K$, such that $K=\bigcup^N_{i=1}K_i$ and $K_i\bigcap K_j=\emptyset$ for $i\neq j$.
    For each element two sets of points are considered: a set located on the boundary, $\partial K$, called the flux points; and a second set called solution points, both such that $\mathbf{x}\in K_i$. The number of solution points is equal to the number of polynomial bases in the approximation space, and the number of flux points is equal to the number of bases in the trace of the approximation space. In one dimension, with an approximation space $\mathbb{P}_k$, there are $k+1$ solution points and $2$ flux points. Lagrange polynomials for the solution and discontinuous flux can then be constructed. To enforce conservation, the discontinuous flux must be made continuous, and in FR the following procedure is used: 
    \begin{equation}\label{eq:fr1d}
        \px{f}{x} \approx J_i^{-1}\left[\px{f^\delta_i}{\xi} + (f^\mathrm{num}_L - f^\delta_L)\dx{h_L}{\xi} + (f^\mathrm{num}_R - f^\delta_R)\dx{h_R}{\xi}\right]
    \end{equation}
    Here, $J^{-1}_i$ is the inverse of the spatial Jacobian. This is used as it is more efficient for interpolation and differentiation operators to work in a reference domain $\hat{K}$ parameterised by $\xi$. Assuming affine elements, we can define the transformation $J_i:\hat{K}\mapsto K$. The last two terms on right-hand side of \cref{eq:fr1d} are the corrections to the flux which ensure conservation. The terms $f^\mathrm{num}_L$ and $f^\mathrm{num}_R$ are common numerical fluxes at the left and right interfaces, respectively, and $f^\delta_L$ and $f^\delta_R$ are the interpolated discontinuous fluxes at the left and right interfaces. Finally, the functions $h_L$ and $h_R$ are the left and right correction functions, with the boundary conditions that they equal one at their respective interfaces, and zero at their opposite interfaces. More detail on the correction functions will be given in the subsequent sub-section. 
    
    Once the continuous gradient of the flux is approximated, the method of lines can be used with an integration method such as explicit Runge--Kutta, or a more complex implicit approach can be used, such as those in \citet{Wang2020}. For a more detailed introduction to the FR method, the works of \citet{DeGrazia2014} and \citet{Abe2015} are recommended, along with the references therein.

\subsection{Correction Functions}
    Since the inception of the FR method~\citep{Huynh2007}, it has been observed that changing the correction function can have a noticeable effect on the scheme's numerical properties. The first continuous set of correction functions was introduced by \citep{Vincent2010}, where a stability proof in one dimension was set out for all functions comprising the set. These functions, parameterised by a single variable, $c$, have the definition:
    \begin{subequations}
        \begin{align}
            h_L &= \frac{(-1)^k}{2}\left(\psi_k - \frac{\eta_k\psi_{k-1} + \psi_{k+1}}{1+\eta_k} \right),\\
            h_R & = \frac{1}{2}\left(\psi_k + \frac{\eta_k\psi_{k-1} + \psi_{k+1}}{1+\eta_k} \right),
        \end{align}
    \end{subequations}
    with the constants:
    \begin{equation}
        \eta_k(c) = \frac{c(2k+1)(a_kk!)}{2}, \quad a_k=\frac{(2k)!}{2^k(k!)^2}, \quad \forall \;c\in\{c\in\mathbb{R};-1<\eta_k(c)<\infty\}.
    \end{equation}
    Here, $\psi_i$ is the $i^\text{th}$ order Legendre polynomial. To construct correction functions for hyper-cube elements such as quadrilaterals and hexahedrons, a tensor product construction of one dimensional functions has typically been used. However, for triangular elements~\citep{Castonguay2011} an analogous proof to that used in 1D was constructed, enabling stable correction functions to be found without a tensor product formulation.
    
    An alternative methodology to define stable correction functions was introduced by \citet{Vincent2015} and later formalised within the summation-by-parts (SBP) framework by \citet{Ranocha2016}. These works only focused on one-dimensional schemes, but they showed the utility of the discrete SBP framework in defining stable schemes. To allow the definition of FR on quadrilaterals to be extended, we now introduce the SBP framework.
    
\subsection{Summation-By-Parts}
    Before defining SBP in higher dimensions, consider the following definitions. Take the domain $K\subset \mathbb{R}^d$, and let $u_i$ be an approximation to the exact function $u$ in element $K_i$. The vector $\mathbf{u}_i$ can then be defined, which is the function $u_i$ evaluated at $N_s$ solution points $\mathbf{x}_i=\{\mathbf{x}_{i,j}\}_{i\leq N_s}$. If we then have the Lagrange polynomials in element $K_i$ such that $l_j(\mathbf{x}_{i,k})=\delta_{jk}$ and $u_i = \sum^{N_s}_{j=1}u_i(\mathbf{x}_{i,j})l_j$, a mass matrix can be defined, with entries:
    \begin{equation}\label{eq:mass_matrix}
        \mathbf{M}_{jk} =  \int_K l_j(\mathbf{x})l_k(\mathbf{x})\mathrm{d}\mathbf{x}.
    \end{equation}
    For cardinal axes $x_1, x_2,\dots$, we can also define the differentiation matrices such that:
    \begin{equation}
        \mathbf{D}_{x_1}\mathbf{u}_i = \sum^{N_s}_{j=1}u_i(\mathbf{x}_{i,j})\dx{l_j}{x_1}, \quad \mathbf{D}_{x_2}\mathbf{u}_i = \sum^{N_s}_{j=1}u_i(\mathbf{x}_{i,j})\dx{l_j}{x_2}, \quad \dots
    \end{equation}
    
    Using these operators we can then define summation-by-parts as a discrete analogy of integration-by-parts, with the following definition:
    \begin{definition}[Generalised Summation-By-Parts]\label{def:sbp}
        Let $u\in C^1(K)$ and $U\in (c^1(K))^d$, such that for some nodal point set $\{\mathbf{x}_i\}_{i\leq N}\subset K$ we have $\mathbf{u}_i = u(\mathbf{x}_i)$ and $\mathbf{U}_i = U(\mathbf{x}_i)$, then a set of operators is said to satisfy the generalised SBP property if:
        \begin{equation}\label{eq:sbp}
            \mathbf{MD} + \mathbf{G}^T\hb{M} = \mathbf{L}^T_{\partial K}\mathbf{W}_{\partial K}\mathbf{N}\hb{L}_{\partial K},
        \end{equation}
        where we have the divergence and gradient operators as:
        \begin{equation}
            \mathbf{DU} = [\mathbf{D}_{x_1}, \mathbf{D}_{x_2}, \dots]\mathbf{U} \approx \nabla\cdot U \quad \mathrm{and} \quad \mathbf{Gu} = 
            \begin{bmatrix}
                \mathbf{D}_{x_1}\\ \mathbf{D}_{x_2}\\\vdots
            \end{bmatrix}\mathbf{u} \approx \nabla u.
        \end{equation}
        Then defining the interpolation $\mathbf{L}_{\partial K}:K\mapsto\partial K$, and boundary mass matrix, $\mathbf{W}_{\partial K}$, such that:
        \begin{equation}
            \mathbf{u}^T_i\mathbf{L}_{\partial K}\mathbf{W}_{\partial K}\mathbf{N}\hb{L}_{\partial K}\mathbf{U}_i = \int_{\partial K}u_iU_i\cdot \mathbf{n}_i\mathrm{d}s,
        \end{equation}
        where $\mathbf{n}$ is a vector function of outwards facing normals at the surface, and $\mathbf{N}$ is a matrix of these normals at the flux points. Here, we use the notation for the Kronecker product with the identity of:
        \begin{equation}
            \hb{B} = \mathbf{B}\otimes\mathbf{I}_d.
        \end{equation}
    \end{definition}
    \begin{remark}
        The definition of the mass matrix given in \cref{eq:mass_matrix} fully integrates the basis, however in many applications a quadrature is used instead of explicitly calculating the mass matrix. From \cref{eq:sbp} it is clear that the mass matrix has to have sufficient accuracy to be able to accurately integrate $\mathbf{u}^T\mathbf{MDu}$, however for some quadratures this is not always possible. Using the works of \citet{Chan2018} and \citet{Trojak2022a}, this problem can be remedied by using a second set of points which do possess sufficient strength. In the context of FR, this additional point set is only required during the operator construction. 
    \end{remark}
    
    With these operators established, the FR method in multiple dimensions can then be rewritten as:
    \begin{equation}
        \nabla\cdot\mathbf{F} \approx \mathbf{D}\mathbf{F} + \mathbf{C}\left((\mathbf{n}\cdot\mathbf{F}^\mathrm{num}) -  \mathbf{N}\hb{L}_\partial\mathbf{F}\right),
    \end{equation}
    where $\mathbf{N}$ is a matrix of outwards facing normals and $\mathbf{C}$ is the correction matrix. This matrix is the discrete analogue of the gradient of the correction function terms in \cref{eq:fr1d}.
    
    In this work we will often work with the modal form of operators. This is due to their relative sparsity compared to the nodal form. Transformation between the modal and nodal representations is performed by the Vandermonde matrix, $\mathbf{V}$, as:
    \begin{equation}
        \mathbf{u} = \mathbf{V}\tb{u},
    \end{equation}
    where $\tb{u}$ is a vector of modal coefficients. An operator matrix, $\mathbf{B}$, is transformed to modal form with:
    \begin{equation}
        \tb{B} = \mathbf{V}^{-1}\mathbf{BV}.
    \end{equation}
    In this work, a tilde is used to denote a matrix or vector in the \emph{modal} representation.
    
\subsection{Polynomial Basis}
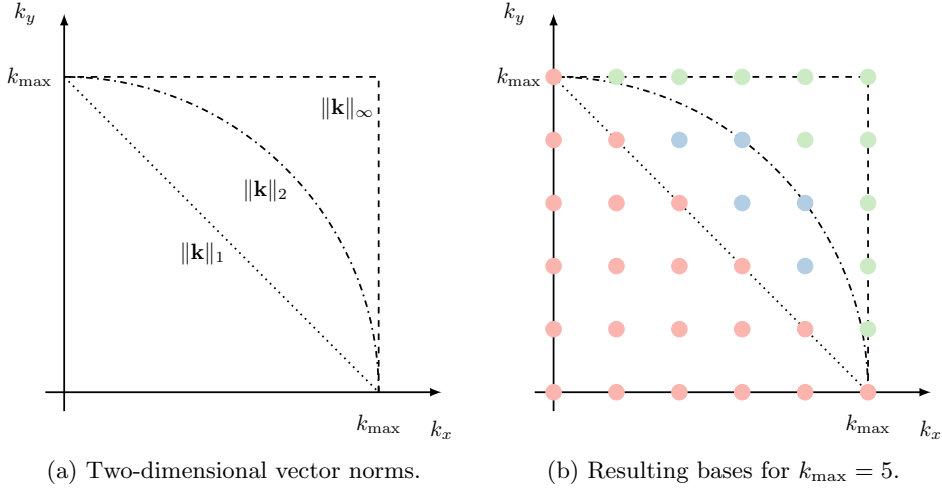
\begin{figure}[tbhp]
        \centering
        \subfloat[Two-dimensional vector norms.]{\adjustbox{width=0.4\linewidth,valign=b}{\input{figs/tikz/order_basis}}}
        ~
        \subfloat[Resulting bases for $k_\mathrm{max}=5$.]{\adjustbox{width=0.4\linewidth,valign=b}{\input{figs/tikz/order_basis_k3}}}
        \caption{\label{fig:basis_order}Diagram of two-dimensional basis orders: total order $\|\mathbf{k}\|_1\leq k_\mathrm{max}$, Euclidean order $\|\mathbf{k}\|_2\leq k_\mathrm{max}$, and maximal order $\|\mathbf{k}\|_\infty\leq k_\mathrm{max}$.}
    \end{figure}

    A systematic way to define a polynomial basis can be achieved through the $L_p$ norm of a vector of orders. This is the method used by \citet{Trefethen2017}, and examples are shown diagrammatically in \cref{fig:basis_order} for two dimensions, where $\mathbf{k}$ is a vector of the basis orders. For example, the basis $\psi_1(x)\psi(y)_2$ would have the vector $[1,2]^T$, where $\psi_i$ is an $i^\text{th}$ order Legendre polynomial. Shown in \cref{fig:basis_order} are the modes required for a total order, Euclidean order, and maximal order basis --- these three bases will form the focus of this work. As outlined in the introduction, on quadrilaterals, maximal order bases have been previously used almost exclusively. One reason for this is that it fits naturally with the element topology. For example, with four corner nodes, the spatial Jacobian can be defined fully in the $k=1$ maximal order basis, \ie bases $1$, $x$, $y$, and $xy$. 
    
    Other basis functions can be chosen --- such as rational functions or radial basis functions. However, except to address some specific deficiencies, these schemes are not widely used due to the additional computational complexities they add, with little benefit in the majority of cases~\citep{Powell1981,Watson2022}.
    

%% file: figs/tikz/order_basis.tex
\begin{tikzpicture}[scale=1]

    \begin{scope}[on behind layer]
        \draw[black, -latex, thick] (0,-0.3) -- (0,6) node [left, text=black, shift={(-0.3, 0.0)}] {\small{$k_y$}};
        \draw[black, -latex, thick] (-0.3,0) -- (6,0) node [below, text=black, shift={( 0.0,-0.3)}] {\small{$k_x$}};
    \end{scope}
    
    \draw[black, -, thick, dashed] (0,5) -- (5,5) -- (5,0);
    
    \draw[black, -, thick, dashdotted] (5,0) arc (0:90:5);
    \draw[black, -, thick, dotted] (0,5) -- (5,0);

    \node[] at (5,-0.5) {\small{$k_\mathrm{max}$}};
    \node[] at (-0.55,5) {\small{$k_\mathrm{max}$}};
    
    \node[] at (2.2,2.2) {\small{$\|\mathbf{k}\|_1$}};
    \node[] at (3.2,3.2) {\small{$\|\mathbf{k}\|_2$}};
    \node[] at (4.5,4.5) {\small{$\|\mathbf{k}\|_\infty$}};
    
\end{tikzpicture}

%% file: figs/tikz/order_basis_k3.tex
\begin{tikzpicture}[scale=1]

    \begin{scope}[on behind layer]
        \draw[black, -latex, thick] (0,-0.3) -- (0,6) node [left, text=black, shift={(-0.3, 0.0)}] {\small{$k_y$}};
        \draw[black, -latex, thick] (-0.3,0) -- (6,0) node [below, text=black, shift={( 0.0,-0.3)}] {\small{$k_x$}};
    \end{scope}
    
    \draw[black, -, thick, dashed] (0,5) -- (5,5) -- (5,0);
    
    \draw[black, -, thick, dashdotted] (5,0) arc (0:90:5);
    \draw[black, -, thick, dotted] (0,5) -- (5,0);

    \node[] at (5,-0.5) {\small{$k_\mathrm{max}$}};
    \node[] at (-0.55,5) {\small{$k_\mathrm{max}$}};
    
    \begin{scope}[on above layer]
        \draw[fill=Pastel1-A, Pastel1-A] (0,0) circle (3.5pt);
        \draw[fill=Pastel1-A, Pastel1-A] (1,0) circle (3.5pt);
        \draw[fill=Pastel1-A, Pastel1-A] (2,0) circle (3.5pt);
        \draw[fill=Pastel1-A, Pastel1-A] (3,0) circle (3.5pt);
        \draw[fill=Pastel1-A, Pastel1-A] (4,0) circle (3.5pt);
        \draw[fill=Pastel1-A, Pastel1-A] (5,0) circle (3.5pt);
        
        \draw[fill=Pastel1-A, Pastel1-A] (0,1) circle (3.5pt);
        \draw[fill=Pastel1-A, Pastel1-A] (1,1) circle (3.5pt);
        \draw[fill=Pastel1-A, Pastel1-A] (2,1) circle (3.5pt);
        \draw[fill=Pastel1-A, Pastel1-A] (3,1) circle (3.5pt);
        \draw[fill=Pastel1-A, Pastel1-A] (4,1) circle (3.5pt);
        
        \draw[fill=Pastel1-A, Pastel1-A] (0,2) circle (3.5pt);
        \draw[fill=Pastel1-A, Pastel1-A] (1,2) circle (3.5pt);
        \draw[fill=Pastel1-A, Pastel1-A] (2,2) circle (3.5pt);
        \draw[fill=Pastel1-A, Pastel1-A] (3,2) circle (3.5pt);
        
        \draw[fill=Pastel1-A, Pastel1-A] (0,3) circle (3.5pt);
        \draw[fill=Pastel1-A, Pastel1-A] (1,3) circle (3.5pt);
        \draw[fill=Pastel1-A, Pastel1-A] (2,3) circle (3.5pt);
        
        \draw[fill=Pastel1-A, Pastel1-A] (0,4) circle (3.5pt);
        \draw[fill=Pastel1-A, Pastel1-A] (1,4) circle (3.5pt);
        
        \draw[fill=Pastel1-A, Pastel1-A] (0,5) circle (3.5pt);
        
        \draw[fill=Pastel1-B, Pastel1-B] (3,3) circle (3.5pt);
        \draw[fill=Pastel1-B, Pastel1-B] (2,4) circle (3.5pt);
        \draw[fill=Pastel1-B, Pastel1-B] (3,4) circle (3.5pt);
        \draw[fill=Pastel1-B, Pastel1-B] (4,2) circle (3.5pt);
        \draw[fill=Pastel1-B, Pastel1-B] (4,3) circle (3.5pt);
        
        \draw[fill=Pastel1-C, Pastel1-C] (5,1) circle (3.5pt);
        \draw[fill=Pastel1-C, Pastel1-C] (5,2) circle (3.5pt);
        \draw[fill=Pastel1-C, Pastel1-C] (5,3) circle (3.5pt);
        \draw[fill=Pastel1-C, Pastel1-C] (5,4) circle (3.5pt);
        \draw[fill=Pastel1-C, Pastel1-C] (5,5) circle (3.5pt);
        
        \draw[fill=Pastel1-C, Pastel1-C] (4,4) circle (3.5pt);
        \draw[fill=Pastel1-C, Pastel1-C] (1,5) circle (3.5pt);
        \draw[fill=Pastel1-C, Pastel1-C] (2,5) circle (3.5pt);
        \draw[fill=Pastel1-C, Pastel1-C] (3,5) circle (3.5pt);
        \draw[fill=Pastel1-C, Pastel1-C] (4,5) circle (3.5pt);
        
    \end{scope}
    
    %\node[] at (2.2,2.2) {\small{$\|\mathbf{k}\|_1$}};
    %\node[] at (3.2,3.2) {\small{$\|\mathbf{k}\|_2$}};
    %\node[] at (4.5,4.5) {\small{$\|\mathbf{k}\|_\infty$}};
    
\end{tikzpicture}

%% file: linear.tex
\section{Linear Stability}\label{sec:linear}
    In the works of \citet{Vincent2015}, \citet{Ranocha2016}, and \citet{Trojak2022a}, the linear stability of flux reconstruction has been explored. The main result of those works is the following lemma for the linear stability of the FR method:
    \begin{lemma}[Linear Stability]\label{lem:lin_stab}
        For flux reconstruction applied to \cref{eq:advect1d} with $\mathbf{f}=\mathbf{F} = \mathbf{a}\otimes\mathbf{u}$, then satisfying the conditions that:
        \begin{subequations}\label{eq:la_cond_main}
            \begin{align}
                \mathbf{Q} &= \mathbf{Q}^T,\label{eq:la_cond4}\\
                (\mathbf{QD}) &= -(\mathbf{QD})^T,\label{eq:la_cond5}\\
                \mathbf{v}^T(\mathbf{M} + \mathbf{Q})\mathbf{v} &> 0,
            \end{align}
        \end{subequations}
        and 
        \begin{equation}\label{eq:la_cond7}
            \mathbf{C} = (\mathbf{M} + \mathbf{Q})^{-1}\mathbf{L}^T_\partial\mathbf{W}_\partial,
        \end{equation}
        with numerical flux such that: 
        \begin{subequations}
            \begin{align}
                (\mathbf{n}\cdot F)^{\mathrm{num} +}_j &= \half(\mathbf{n}_j^+\cdot\mathbf{a})(u_j^+ + u_j^-) - \half\kappa|\mathbf{n}_j^+\cdot\mathbf{a}|(u_j^- - u_j^+), \quad \text{and}\\
                (\mathbf{n}\cdot F)^{\mathrm{num} -}_j &= \half(\mathbf{n}_j^-\cdot\mathbf{a})(u_j^- + u_j^+) - \half\kappa|\mathbf{n}_j^-\cdot\mathbf{a}|(u_j^+ - u_j^-), \quad \text{for} \quad \kappa \in [0,1],
            \end{align}
        \end{subequations}
        means the scheme is linearly stable, in that:
        \begin{equation}
            \dx{}{t}\|\mathbf{u}\|_{M+Q}^2 \leq 0.
        \end{equation}
    \end{lemma}
    \begin{proof}
        For a proof see \citet{Trojak2022a}.
    \end{proof}
    
    The conditions set out in \cref{eq:la_cond_main,eq:la_cond7} allow for a parameterised $\mathbf{Q}$ that defines a continuous set of stable FR schemes to be found. The reduction of a generic $\mathbf{Q}$ matrix to enforce these conditions can be performed in a symbolic manipulation toolbox, and, by doing so, a general framework can be produced to find stable sets of FR schemes.
    
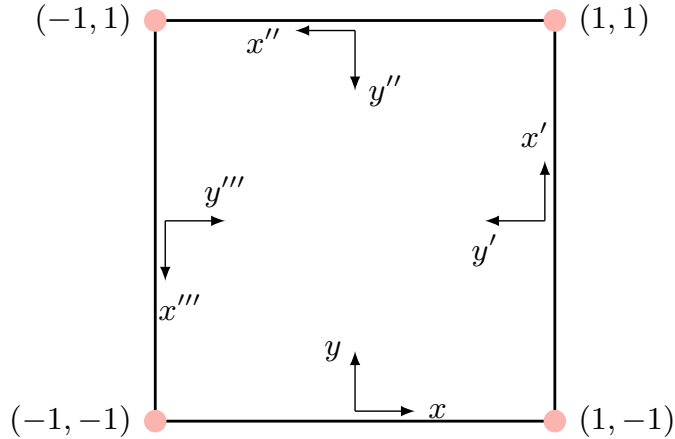
\begin{figure}[tbhp]
        \centering
        \adjustbox{width=0.6\linewidth,valign=b}{\input{figs/tikz/ref_square}}
        \caption{\label{fig:ref_square}Reference quadrilateral and the four face-relative coordinate systems.}
    \end{figure}

    In addition to these conditions, it is assumed that the numerical properties of the method should be independent of the node ordering. Therefore, additional symmetry conditions are required for $\mathbf{Q}$ such that, for the four reference axes shown in \cref{fig:ref_square}, $\mathbf{Q}$ is independent of a particular frame of reference. 
    
    To achieve the desired symmetry properties, we first start by defining a transformation matrix from one reference frame to another, $\mathbf{T}$, and then enforce the following condition:
    \begin{equation}
        \tb{T}_{ab}\tb{Q} = \tb{Q}\tb{T}_{ab},
    \end{equation}
    here enforced in the modal representation. The matrix $\mathbf{T}_{ab}$ transforms a vector from reference frame $a$ to frame $b$. In later sections, we will go on to explore alternative bases, for which rotationally symmetric point layouts are not possible. In these situations, a certain degree of anisotropy will have to be accepted, and at least with these symmetry conditions, the methods will be as symmetric as possible. Care should be taken when enforcing the symmetry conditions to not over-constrain $\mathbf{Q}$. For a quadrilateral, this means that only two rotations need to be enforced, as the remaining rotational and axial symmetries can be expressed in terms of just two rotations.

%% file: figs/tikz/ref_square.tex
\begin{tikzpicture}[scale=2]
    
    \begin{scope}
        \draw[fill={Pastel1-A}, color={Pastel1-A}] (-1,-1) circle (1.5pt) node [left, text=black, shift={(-0.1,0.0)}] {\footnotesize{$\left(-1,-1\right)$}};
        
        \draw[fill={Pastel1-A}, color={Pastel1-A}] ( 1,-1) circle (1.5pt) node [right, text=black, shift={(0.1,0.0)}] {\footnotesize{$\left( 1,-1\right)$}};
        
        \draw[fill={Pastel1-A}, color={Pastel1-A}] ( 1, 1) circle (1.5pt) node [right, text=black, shift={(0.1,0.0)}] {\footnotesize{$\left( 1, 1\right)$}};
        
        \draw[fill={Pastel1-A}, color={Pastel1-A}] (-1, 1) circle (1.5pt) node [left, text=black, shift={(-0.1,0.0)}] {\footnotesize{$\left(-1, 1\right)$}};
        
    \end{scope}
    
    \begin{scope}[on behind layer]
        \draw[thick, black] (-1,-1) -- (1,-1) -- (1,1) -- (-1,1) -- cycle;
    \end{scope}
    
    \draw[black, -latex] (0,-0.95) -- (0,-0.65) node[left] {\footnotesize{$y$}};
    \draw[black, -latex] (0,-0.95) -- (0.3,-0.95) node[right] {\footnotesize{$x$}};
    
    \draw[black, -latex] (0.95,0) -- (0.65,0) node[below, shift={(0.0,0)}] {\footnotesize{$y^{\prime}$}};
    \draw[black, -latex] (0.95,0) -- (0.95,0.3) node[above, shift={(-0.1,0)}] {\footnotesize{$x^{\prime}$}};
    
    \draw[black, -latex] (0,0.95) -- (0,0.65) node[right, shift={(0.0,0)}] {\footnotesize{$y^{\prime\prime}$}};
    \draw[black, -latex] (0,0.95) -- (-0.3,0.95) node[left, shift={(-0,-0.1)}] {\footnotesize{$x^{\prime\prime}$}};
    
    \draw[black, -latex] (-0.95,0) -- (-0.65,0) node[above, shift={(0.0,0)}] {\footnotesize{$y^{\prime\prime\prime}$}};
    \draw[black, -latex] (-0.95,0) -- (-0.95,-0.3) node[below, shift={(0.15,0)}] {\footnotesize{$x^{\prime\prime\prime}$}};
    
\end{tikzpicture}

%% file: erfr_quad_mo.tex
\section{Extended-range FR for quadrilaterals}\label{sec:quad_mo}
    The overwhelming majority of polynomial finite element methods when applied to quadrilaterals use a maximal order basis, \ie $\|\mathbf{k}\|_\infty\leq k_\mathrm{max}$. To define an extended range of stable FR method in this case, the techniques of \cref{sec:linear} can be applied. There are many possible options for the point sets. It has been shown that a tensor product of Gauss--Lobatto points is Fekete optimal~\citep{Bos2001}, and that a tensor product of Chebyshev points is near optimal in a Lebesgue sense~\citep{Caliari2005}. However, for methods such as FR, it is has been shown in one dimension that Gauss--Legendre points are optimal, and it has been suggested that this extends to higher dimensions via a tensor product~\citep{Witherden2021}. 
    
    The reference element for the quadrilateral used in this work is shown in \cref{fig:ref_square}, and the maximal order orthogonal basis is organised as:
    \begin{equation}
        \phi_i(x,y) = \psi_v(x)\psi_w(y), \quad \mathrm{for} \quad i=w(k+1)+v+1, \quad \mathrm{and} \quad 0\leq v,w \leq k.
    \end{equation}
    
\subsection{$k=2$}
    Starting at $k=2$, the conditions set out in \cref{lem:lin_stab} and the symmetry conditions can be enforced on a matrix, to find that applicable $\tb{Q}$ matrices have the form:
    \begin{equation}
        \tb{Q} = \left[ \begin{array}{cccccccc}
            \mathbf{0}  & \multicolumn{7}{c}{\mathbf{0}} \\
            \multirow{7}{*}{$\boldsymbol{0}$}
            & 0 & 0 & 0 & 0 & 0 & 0 & -3q_1 \\
            & 0 & 0 & 0 & 0 & 0 & 0 & 0 \\
            & 0 & 0 & 0 & 0 & 0 & 0 & 0 \\
            & 0 & 0 & 0 & q_1 & 0 & 0 & 0 \\
            & 0 & 0 & 0 & 0 & 0 & 0 & -3q_1 \\
            & 0 & 0 & 0 & 0 & 0 & q_1 & 0 \\
            & -3q_1 & 0 & 0 & 0 & -3q_1 & 0 & q_0
        \end{array}\right].
    \end{equation}
    For stability, it is required that $\tb{M}+\tb{Q}$ is positive definite in order to induce a valid norm. Therefore, this imposes some conditions on the values of $\tb{Q}$. These can be straightforwardly found via the Cholesky factorisation, and for $k=2$ the conditions are:
    \begin{equation}
            q_1 > -4/15 \quad \mathrm{and} \quad 50q_0 -1125q_1^2 + 8 > 0.
    \end{equation}
    
\subsection{$k=3$}
    The analysis can be repeated for $k=3$, to obtain:
    \begin{equation}
        \tb{Q} = \left[ \begin{array}{cccccccccccc}
            \mathbf{0}  & \multicolumn{11}{c}{\mathbf{0}} \\
            \multirow{11}{*}{$\boldsymbol{0}$} & 0 & 0 & 0 & 0 & 0 & 0 & 0 & 0 & 0 & 0 & q_2 \\
             & 0 & 0 & 0 & 0 & 0 & 0 & 0 & 0 & 0 & -3q_2/5 & 0 \\
             & 0 & 0 & 0 & 0 & 0 & 0 & 0 & 0 & q_2 & 0 & -5q_1/3 \\
             & 0 & 0 & 0 & 0 & 0 & 0 & 0 & 0 & 0 & 0 & 0 \\
             & 0 & 0 & 0 & 0 & 0 & 0 & -3q_2/5 & 0 & 0 & 0 & 0 \\
             & 0 & 0 & 0 & 0 & 0 & 9q_2/25 & 0 & 0 & 0 & 0 & 0 \\
             & 0 & 0 & 0 & 0 & -3q_2/5 & 0 & q_1 & 0 & 0 & 0 & 0 \\
             & 0 & 0 & 0 & 0 & 0 & 0 & 0 & 0 & 0 & 0 & 0 \\
             & 0 & 0 & q_2 & 0 & 0 & 0 & 0 & 0 & 0 & 0 & -5q_1/3 \\
             & 0 & -3q_2/5 & 0 & 0 & 0 & 0 & 0 & 0 & 0 & q_1 & 0 \\
             & q_2 & 0 & -5q_1/3 & 0 & 0 & 0 & 0 & 0 & -5q_1/3 & 0 & q_0
        \end{array}\right]
    \end{equation}
    for the conditions on stability that:
    \begin{subequations}
        \begin{align}
            q_2^2 &< \frac{16}{441},\\
            - 189q_2^2 + 140q_1 &> -16,\\
            (4 - 21q_2)(1008q_2 + 2352q_0 + 12348q_2q_0 - 5292q_2^2 - 27783q_2^3 - 68600q_1^2 + 192) &>0.
        \end{align}
    \end{subequations}
    
    This procedure can be continued for any order, $k$, to recover the $\hb{Q}$ matrix and stability conditions. The results for $k \geq 4$ are cumbersome and are therefore excluded for brevity and typesetting constraints.

\subsection{Tensor-product schemes}
    In the earlier works on the topic of stable FR schemes for quadrilaterals, correction functions were constructed using a tensor product of stable one-dimensional schemes. We wish to understand if these tensor-product constructions can be found as a subset of the schemes defined here. 
    
    To do this we first consider the modal presentation of the one-dimensional class of \citet{Vincent2010}, which can be used to formulate a tensor-product modal correction matrix. For the case of $k=2$ this leads to the $\tb{C}$ matrix:
    \begin{equation}
        \tb{C}_{tp} = \left[\begin{array}{cccccccccccc}
            1/2 & 0 & 0 & 1/2 & 0 & 0 & 1/2 & 0 & 0 & 1/2 & 0 & 0 \\
            0 & 1/2 & 0 & 3/2 & 0 & 0 & 0 & -1/2 & 0 & -3/2 & 0 & 0 \\
            0 & 0 & 1/2 & \theta & 0 & 0 & 0 & 0 & 1/2 & \theta & 0 & 0 \\
            -3/2 & 0 & 0 & 0 & 1/2 & 0 & 3/2 & 0 & 0 & 0 & -1/2 & 0 \\
            0 & -3/2 & 0 & 0 & 3/2 & 0 & 0 & -3/2 & 0 & 0 & 3/2 & 0 \\
            0 & 0 & -3/2 & 0 & \theta & 0 & 0 & 0 & 3/2 & 0 & -\theta & 0 \\
            \theta & 0 & 0 & 0 & 0 & 1/2 & \theta & 0 & 0 & 0 & 0 & 1/2 \\
            0 & \theta & 0 & 0 & 0 & 3/2 & 0 & -\theta & 0 & 0 & 0 & -3/2 \\
            0 & 0 & \theta & 0 & 0 & \theta & 0 & 0 & \theta & 0 & 0 & \theta
        \end{array}\right],
    \end{equation}
    for $\theta = 5/(45c + 2)$. Attempts can then be made to solve the following system to find a valid $\tb{Q}$:
    \begin{equation}\label{eq:tp_system}
        \tb{Q}\tb{C}_{tp} = -\tb{M}(\tb{C}_{tp} - \tb{C}_{DG}),
    \end{equation}
    where $\tb{C}_{DG}$ is the DG correction matrix, found from $\mathbf{C}_{DG} = \mathbf{M}^{-1}\mathbf{L}^T_\partial\mathbf{W}_\partial$. This substitution is used in \cref{eq:la_cond7} as it gives a simpler system to solve. Looking for solutions, only one is found: when $c=0$ and $\tb{Q}=0$. 
    
    Repeating this for analysis for the extended range of stable 1D FR schemes presented by \citet{Vincent2015}, we find the tensor-product modal correction matrix for $k=2$ as:
    \begin{equation}
        \tb{C}_{tp} = \left[\begin{array}{cccccccccccc}
    -1/2 & 0 & 0 & 1/2 & 0 & 0 & 1/2 & 0 & 0 & -1/2 & 0 & 0 \\
     0 & -1/2 & 0 & -\theta_0 & 0 & 0 & 0 & -1/2 & 0 & -\theta_0 & 0 & 0 \\
     0 & 0 & -1/2 & \theta_1 & 0 & 0 & 0 & 0 & 1/2 & -\theta_1 & 0 & 0 \\
     -\theta_0 & 0 & 0 & 0 & 1/2 & 0 & -\theta_0 & 0 & 0 & 0 & 1/2 & 0 \\
     0 & -\theta_0 & 0 & 0 & -\theta_0 & 0 & 0 & \theta_0 & 0 & 0 & \theta_0 & 0 \\
     0 & 0 & -\theta_0 & 0 & \theta_1 & 0 & 0 & 0 & -\theta_0 & 0 & \theta_1 & 0 \\
     -\theta_1 & 0 & 0 & 0 & 0 & 1/2 & \theta_1 & 0 & 0 & 0 & 0 & -1/2 \\
     0 & -\theta_1 & 0 & 0 & 0 & -\theta_0 & 0 & -\theta_1 & 0 & 0 & 0 & -\theta_0 \\
     0 & 0 & -\theta_1 & 0 & 0 & \theta_1 & 0 & 0 & \theta_1 & 0 & 0 & -\theta_1
        \end{array}\right],
    \end{equation}
    with
    \begin{equation}
        \theta_0 = (63c_0 + 105c_1 + 18)/\Psi, \quad \Psi = 175c_1^2 - 42c_0 - 12, \quad \text{and} \theta_1 = 5/(5c_1 + 2).
    \end{equation}
    Once more, solutions to the system shown in \cref{eq:tp_system} can be sought, whereupon it is found that no solutions exist except for $c_0=c_1=0$ --- the DG solution.
    This leads us to the following proposition: for quadrilateral elements, a correction matrix that is a tensor-product of a one-dimensional correction function is not a form of linearly stable filtered DG scheme --- with the exception of DG itself --- although norms can exist where monotonic decay is observed.

%% file: erfr_quad_to.tex
\section{Total order basis}\label{sec:quad_to}
    Rather than the typical maximal order basis, if instead a total order basis is used, such that $\|\mathbf{k}\|_1\leq k_\mathrm{max}$, then a new set of stable FR schemes can be recovered. This basis is analogous to that used on triangular elements. A key requirement for finite element numerical methods is that the approximation space on the element boundary is the trace of the approximation space of the element. An advantage of hyper-cube topologies, such as the quadrilateral, is that it is trivial to show that for $\|\mathbf{k}\|_p\leq k_\mathrm{max}$ this is true for $0<p\leq\infty$.
    
    \begin{table}[tbhp]
        \centering
        \begin{tabular}{r r r r}
            \toprule
            \multirow{2}{*}{$k_\mathrm{max}$} & \multicolumn{3}{c}{$n_b$} \\ \cmidrule(l){2-4}
            & $\|\mathbf{k}\|_1$ & $\|\mathbf{k}\|_2$ & $\|\mathbf{k}\|_\infty$ \\\midrule
            1 & 3 & 3 & 4 \\
            2 & 6 & 6 & 9 \\
            3 & 10 & 11 & 16 \\
            4 & 15 & 17 & 25 \\
            5 & 21 & 26 & 36 \\
            6 & 28 & 35 & 49 \\
            \bottomrule
        \end{tabular}
        \caption{\label{tab:basis_points}Number of bases, $n_b$, for different norms in two dimensions.}
    \end{table}
    
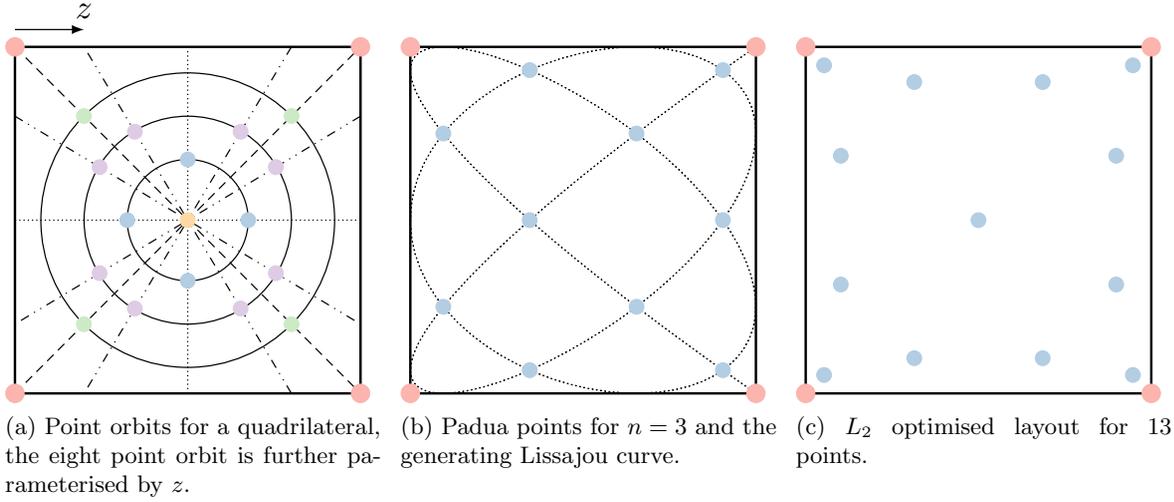
\begin{figure}[tbhp]
        \centering
        \subfloat[Point orbits for a quadrilateral, the eight point orbit is further parameterised by $z$.]{\label{fig:quad_orbits}\adjustbox{width=0.32\linewidth,valign=b}{\input{figs/tikz/quad_orbits}}}
        ~
        \subfloat[Padua points for $n=3$ and the generating Lissajou curve.]{\label{fig:padua}\adjustbox{width=0.32\linewidth,valign=b}{\input{figs/tikz/quad_padua}}}
        ~
        \subfloat[$L_2$ optimised layout for 13 points.]{\label{fig:euclid}\adjustbox{width=0.32\linewidth,valign=b}{\input{figs/tikz/quad_euclid}}}
        \caption{\label{fig:quad_points} Various solution point layouts on the reference quadrilateral.}
    \end{figure}

    In previous literature it has often been taken as axiomatic that the solution points should be placed such that their location is independent of the corner-node ordering. In the work of \citet{Witherden2015} quadratures were found by enforcing this symmetry through orbit groups. For a square there are four such groups, and these groups are shown diagrammatically in \cref{fig:quad_orbits}. However, from the number of basis functions for a given $k_\mathrm{max}$ shown in \cref{tab:basis_points}, it is apparent that the number of total order bases can not always be recovered using these orbits. For an example, consider $k_\mathrm{max}=2$ --- with six bases, the closest symmetric point layout would have five points.
    
    One alternative to symmetric point layouts for total order are the Padua points~\citep{Caliari2005}, an example of which are shown in \cref{fig:padua}. These points have several attractive properties: they provably minimise the Lebesgue measure on the domain $[-1,1]^2$; are unisolvent for arbitrary orders; and have $(n+1)(n+2)/2$ points, \ie they have the same number of points as a total order basis. However, they lack the full rotational symmetry of \cref{fig:quad_orbits}.
    
    With the total order basis introduced, we now enumerate some of the set of linearly stable FR methods on quadrilaterals with a total order basis. 
    
\subsection{$k=2$}
    Starting with $k=2$, enforcing the conditions on stability as presented in \cref{lem:lin_stab}, we find that $\tb{Q}$ can have the form: 
    \begin{equation}
        \tb{Q} = \left[ \begin{array}{cccccc}
            0 & 0 & 0 & 0 & 0 & 0 \\
            0 & 0 & 0 & 0 & 0 & 0 \\
            0 & 0 & q_0 & 0 & 0 & q_2 \\
            0 & 0 & 0 & 0 & 0 & 0 \\
            0 & 0 & 0 & 0 & q_1 & 0 \\
            0 & 0 & q_2 & 0 & 0 & q_0
        \end{array}\right]\quad \mathrm{for} \quad 
        \mathbf{\Phi}_2(x,y) = \begin{bmatrix}
            1 \\
            \psi_0(x)\psi_1(y) \\
            \psi_0(x)\psi_2(y) \\
            \psi_1(x)\psi_0(y) \\
            \psi_1(x)\psi_1(y) \\
            \psi_2(x)\psi_0(y)
        \end{bmatrix}
    \end{equation}
    The condition of $\mathbf{M}+\mathbf{Q}$ being positive definite then leads to the conditions on stability that:
    \begin{equation}
        q_1 > -4/9,\quad q_0 > -4/9, \quad \mathrm{and} \quad (5q_0 + 4)^2 - 25q_2^2 >0.
    \end{equation}
    These conditions can be straightforwardly recovered from the condition that the Cholesky factorisation of a positive definite matrix has positive-real values on the leading diagonal.
    
\subsection{$k=3$}
    Repeating this process for $k=3$:
    \begin{equation}
        \tb{Q} = \left[ \begin{array}{cccccccc}
            \mathbf{0}  & \multicolumn{7}{c}{\mathbf{0}} \\
            \multirow{7}{*}{$\boldsymbol{0}$} 
            & q_0 & 0 & 0 & 0 & 0 & q_2 & 0 \\
            & 0 & 0 & 0 & 0 & 0 & 0 & 0 \\
            & 0 & 0 & 0 & 0 & 0 & 0 & 0 \\
            & 0 & 0 & 0 & q_1 & 0 & 0 & q_2 \\
            & 0 & 0 & 0 & 0 & 0 & 0 & 0 \\
            & q_2 & 0 & 0 & 0 & 0 & q_1 & 0 \\
            & 0 & 0 & 0 & q_2 & 0 & 0 & q_0
        \end{array}\right] \quad \mathrm{for} \quad \mathbf{\Phi}_3(x,y) = \begin{bmatrix}
             1 \\
             \vdots \\
             \psi_0(x)\psi_3(y) \\
             \psi_1(x)\psi_0(y) \\
             \psi_1(x)\psi_1(y) \\
             \psi_1(x)\psi_2(y) \\
             \psi_2(x)\psi_0(y) \\
             \psi_2(x)\psi_1(y) \\
             \psi_3(x)\psi_0(y)
        \end{bmatrix}
    \end{equation}
    subject to the conditions that:
    \begin{equation}
        q_0 > -4/7,\quad q_1 > -4/15,\quad \mathrm{and}\quad 28q_0 + 105q_0q_1 + 60q_1 - 105q_2^2 + 16 > 0.
    \end{equation}
    As an example of the resulting correction field, \cref{fig:quad_to_k3_dg} shows the $k=3$ correction field for DG FR for two different flux points.
    
    \begin{figure}[tbhp]
        \centering
        %\subfloat[First flux point $(-1,\sqrt{(15 + 2\sqrt{30})35}/35)$.]{\adjustbox{width=0.39\linewidth,valign=b}{\input{figs/tikz/quad_p3_corr_f1_dg_to}}}
        \subfloat[First flux point $(-1,\frac{\sqrt{15+2\sqrt{30}}}{35})$.]{\adjustbox{width=0.39\linewidth,valign=b}{\input{figs/tikz/quad_p3_corr_f1_dg_to}}}
        ~
        %\subfloat[Second flux point $(-1,\sqrt{(15 - 2\sqrt{30})35}/35)$.]{\adjustbox{width=0.56\linewidth,valign=b}{\input{figs/tikz/quad_p3_corr_f2_dg_to}}}
        \subfloat[Second flux point $(-1,\frac{\sqrt{15-2\sqrt{30}}}{35})$.]{\adjustbox{width=0.56\linewidth,valign=b}{\input{figs/tikz/quad_p3_corr_f2_dg_to}}}        
        \caption{\label{fig:quad_to_k3_dg}Divergence of DG correction field for $k=3$ FR on a quadrilateral with total order basis for two flux points, shown in red.}
    \end{figure}
    
\subsection{$k=4$}
    By repeating the process again, the $\tb{Q}$ matrix and stability conditions have been found for $k=4$:
    \begin{equation}
        \tb{Q} = \left[ \begin{array}{cccccccccccc}
            \mathbf{0} & \multicolumn{11}{c}{\mathbf{0}} \\
            \multirow{11}{*}{$\boldsymbol{0}$} 
      & q_0 & 0 & 0 & 0 & 0 & 0 & 0 & q_3 & 0 & 0 & q_5 \\
      & 0 & 0 & 0 & 0 & 0 & 0 & 0 & 0 & 0 & 0 & 0 \\
      & 0 & 0 & 0 & 0 & 0 & 0 & 0 & 0 & 0 & 0 & 0 \\
      & 0 & 0 & 0 & 0 & 0 & 0 & 0 & 0 & 0 & 0 & 0 \\
      & 0 & 0 & 0 & 0 & q_1 & 0 & 0 & 0 & 0 & q_4 & 0 \\
      & 0 & 0 & 0 & 0 & 0 & 0 & 0 & 0 & 0 & 0 & 0 \\
      & 0 & 0 & 0 & 0 & 0 & 0 & 0 & 0 & 0 & 0 & 0 \\
      & q_3 & 0 & 0 & 0 & 0 & 0 & 0 & q_2 & 0 & 0 & q_3 \\
      & 0 & 0 & 0 & 0 & 0 & 0 & 0 & 0 & 0 & 0 & 0 \\
      & 0 & 0 & 0 & 0 & q_4 & 0 & 0 & 0 & 0 & q_1 & 0 \\
      & q_5 & 0 & 0 & 0 & 0 & 0 & 0 & q_3 & 0 & 0 & q_0
        \end{array}\right] \quad \mathrm{for} \quad
        \mathbf{\Phi}_4(x,y) = \begin{bmatrix}
            1 \\
            \vdots\\
     \psi_0(x)\psi_4(y)\\
     \psi_1(x)\psi_0(y)\\
     \psi_1(x)\psi_1(y)\\
     \psi_1(x)\psi_2(y)\\
     \psi_1(x)\psi_3(y)\\
     \psi_2(x)\psi_0(y)\\
     \psi_2(x)\psi_1(y)\\
     \psi_2(x)\psi_2(y)\\
     \psi_3(x)\psi_0(y)\\
     \psi_3(x)\psi_1(y)\\
     \psi_4(x)\psi_0(y)
        \end{bmatrix}
    \end{equation}
    subject to the constraints:
    \begin{subequations}
        \begin{align}
            q_0 &> -4/9,\\
            q_1 &> -4/21,\\
            36q_0 + 225q_0q_2 + 100q_2 - 225q_3^2 + 16 &> 0,\\
            (21q_1+4)^2 - 441q_4^2 &> 0,\\
            (9q_0 - 9q_5 + 4)\left[9q_0(25q_2 + 4) + 25q_2(9 q_5 + 4) + 2(-225q_3^2 + 18 q_5 + 8)\right] &>0.
        \end{align}
    \end{subequations}

%% file: figs/tikz/quad_orbits.tex
\begin{tikzpicture}[scale=2]
    
    \begin{scope}
        \draw[fill={Pastel1-A}, color={Pastel1-A}] (-1,-1) circle (1.5pt);
        
        \draw[fill={Pastel1-A}, color={Pastel1-A}] ( 1,-1) circle (1.5pt);
        
        \draw[fill={Pastel1-A}, color={Pastel1-A}] ( 1, 1) circle (1.5pt);
        
        \draw[fill={Pastel1-A}, color={Pastel1-A}] (-1, 1) circle (1.5pt);
        
    \end{scope}
    
    \begin{scope}[on behind layer]
        \draw[thick, black] (-1,-1) -- (1,-1) -- (1,1) -- (-1,1) -- cycle;
    \end{scope}
    
    \begin{scope}[on behind layer]
        \draw[black, densely dashed] (-1,-1) -- (1,1);
        \draw[black, densely dashed] (-1, 1) -- (1,-1);
        \draw[black] (0,0) circle (0.85);
        
        \draw[black, densely dotted] (0,-1) -- (0,1);
        \draw[black, densely dotted] (-1,0) -- (1,0);
        \draw[black] (0,0) circle (0.35);
        
        \draw[black, dashdotdotted] (0.6,-1) -- (-0.6,1);
        \draw[black, dashdotdotted] (0.6,1) -- (-0.6,-1);
        \draw[black, dashdotdotted] (-1,0.6) -- (1,-0.6);
        \draw[black, dashdotdotted] (1,0.6) -- (-1,-0.6);
        \draw[black] (0,0) circle (0.6);
    \end{scope}
    
    \begin{scope}[on above layer]
        \draw[fill=Pastel1-B,Pastel1-B] (0, 0.35) circle (1.2pt);
        \draw[fill=Pastel1-B,Pastel1-B] (0,-0.35) circle (1.2pt);
        \draw[fill=Pastel1-B,Pastel1-B] ( 0.35,0) circle (1.2pt);
        \draw[fill=Pastel1-B,Pastel1-B] (-0.35,0) circle (1.2pt);
        
        \draw[fill=Pastel1-C,Pastel1-C] ( 0.601, 0.601) circle (1.2pt);
        \draw[fill=Pastel1-C,Pastel1-C] ( 0.601,-0.601) circle (1.2pt);
        \draw[fill=Pastel1-C,Pastel1-C] (-0.601, 0.601) circle (1.2pt);
        \draw[fill=Pastel1-C,Pastel1-C] (-0.601,-0.601) circle (1.2pt);
        
        \draw[fill=Pastel1-D,Pastel1-D] (0.306,0.51) circle (1.2pt);
        \draw[fill=Pastel1-D,Pastel1-D] (0.51,0.306) circle (1.2pt);
        \draw[fill=Pastel1-D,Pastel1-D] (-0.306,0.51) circle (1.2pt);
        \draw[fill=Pastel1-D,Pastel1-D] (-0.51,0.306) circle (1.2pt);
        \draw[fill=Pastel1-D,Pastel1-D] (0.306,-0.51) circle (1.2pt);
        \draw[fill=Pastel1-D,Pastel1-D] (0.51,-0.306) circle (1.2pt);
        \draw[fill=Pastel1-D,Pastel1-D] (-0.306,-0.51) circle (1.2pt);
        \draw[fill=Pastel1-D,Pastel1-D] (-0.51,-0.306) circle (1.2pt);
        
        \draw[fill=Pastel1-E,Pastel1-E] (0,0) circle (1.2pt);
        
        \draw[black,-latex] (-1,1.1) -- (-0.6,1.1) node[left, above] {$z$};
        
    \end{scope}
    
\end{tikzpicture}

%% file: figs/tikz/quad_padua.tex
\begin{tikzpicture}[scale=2]
    
    \begin{scope}
        \draw[fill={Pastel1-A}, color={Pastel1-A}] (-1,-1) circle (1.5pt);
        
        \draw[fill={Pastel1-A}, color={Pastel1-A}] ( 1,-1) circle (1.5pt);
        
        \draw[fill={Pastel1-A}, color={Pastel1-A}] ( 1, 1) circle (1.5pt);
        
        \draw[fill={Pastel1-A}, color={Pastel1-A}] (-1, 1) circle (1.5pt);
        
    \end{scope}
    
    \begin{scope}[on behind layer]
        \draw[thick, black] (-1,-1) -- (1,-1) -- (1,1) -- (-1,1) -- cycle;
        
        \draw[domain=0:360, samples=360, variable=\x, black, densely dotted] plot ({cos(6*\x)}, {cos(5*\x)});
    \end{scope}
    
    \draw[fill=Pastel1-B,Pastel1-B] (-0.80902,-0.5) circle (1.2pt);
    \draw[fill=Pastel1-B,Pastel1-B] (-0.80902,0.5) circle (1.2pt);
    \draw[fill=Pastel1-B,Pastel1-B] (-0.30902,-0.86603) circle (1.2pt);
    \draw[fill=Pastel1-B,Pastel1-B] (-0.30902,6.1232e-17) circle (1.2pt);
    \draw[fill=Pastel1-B,Pastel1-B] (-0.30902,0.86603) circle (1.2pt);
    \draw[fill=Pastel1-B,Pastel1-B] (0.30902,-0.5) circle (1.2pt);
    \draw[fill=Pastel1-B,Pastel1-B] (0.30902,0.5) circle (1.2pt);
    \draw[fill=Pastel1-B,Pastel1-B] (0.80902,-0.86603) circle (1.2pt);
    \draw[fill=Pastel1-B,Pastel1-B] (0.80902,6.1232e-17) circle (1.2pt);
    \draw[fill=Pastel1-B,Pastel1-B] (0.80902,0.86603) circle (1.2pt);
    
\end{tikzpicture}

%% file: figs/tikz/quad_euclid.tex
\begin{tikzpicture}[scale=2]
    
    \begin{scope}
        \draw[fill={Pastel1-A}, color={Pastel1-A}] (-1,-1) circle (1.5pt);
        
        \draw[fill={Pastel1-A}, color={Pastel1-A}] ( 1,-1) circle (1.5pt);
        
        \draw[fill={Pastel1-A}, color={Pastel1-A}] ( 1, 1) circle (1.5pt);
        
        \draw[fill={Pastel1-A}, color={Pastel1-A}] (-1, 1) circle (1.5pt);
        
    \end{scope}
    
    \begin{scope}[on behind layer]
        \draw[thick, black] (-1,-1) -- (1,-1) -- (1,1) -- (-1,1) -- cycle;
    \end{scope}
    
    \draw[fill=Pastel1-B,Pastel1-B] (0,0) circle (1.2pt);
    \draw[fill=Pastel1-B,Pastel1-B] (0.89367,-0.89367) circle (1.2pt);
    \draw[fill=Pastel1-B,Pastel1-B] (0.89367,0.89367) circle (1.2pt);
    \draw[fill=Pastel1-B,Pastel1-B] (-0.89367,0.89367) circle (1.2pt);
    \draw[fill=Pastel1-B,Pastel1-B] (-0.89367,-0.89367) circle (1.2pt);
    \draw[fill=Pastel1-B,Pastel1-B] (0.37165,-0.79694) circle (1.2pt);
    \draw[fill=Pastel1-B,Pastel1-B] (0.79694,0.37165) circle (1.2pt);
    \draw[fill=Pastel1-B,Pastel1-B] (-0.37165,0.79694) circle (1.2pt);
    \draw[fill=Pastel1-B,Pastel1-B] (-0.79694,-0.37165) circle (1.2pt);
    \draw[fill=Pastel1-B,Pastel1-B] (0.79694,-0.37165) circle (1.2pt);
    \draw[fill=Pastel1-B,Pastel1-B] (-0.37165,-0.79694) circle (1.2pt);
    \draw[fill=Pastel1-B,Pastel1-B] (-0.79694,0.37165) circle (1.2pt);
    \draw[fill=Pastel1-B,Pastel1-B] (0.37165,0.79694) circle (1.2pt);
    
\end{tikzpicture}

%% file: figs/tikz/quad_p3_corr_f1_dg_to.tex
    \begin{tikzpicture}[scale=2]
		\begin{axis}[name=plot1,xlabel={$x$},ylabel={$y$},
    		axis line style={latex-latex},
            axis y line=middle,
            axis x line=middle,
            xmode=linear, % not log
            ymode=linear, % not log
		    xtick={0},
		    ytick={0},
		    xticklabels={ },
		    yticklabels={ },
		    ylabel style={rotate=0, shift={(-0.09,0.06)}},
		    xlabel style={rotate=0, shift={(0.03,0.01)}},
    		xmin=-1.1,xmax=1.1,
    		ymin=-1.1,ymax=1.1,
    		width=5cm,
    		height=5cm,
			scale only axis=true,
    		style={font=\normalsize},
    		axis on top,
    % 		colorbar,
    % 		colormap/Blues,
    % 		colorbar style={
    % 		    point meta max=5, point meta min=-2,
    %             ylabel=$\nabla\cdot\mathbf{h}$,
    %             ylabel style={font=\normalsize},
    %             ytick={-4,-2,0,2},
    %             yticklabel style={
    %                 text width=2em,
    %                 align=left,
    %                 font=\normalsize,
    %             }
    %         }
        ]
            \addplot[thick] graphics[xmin=-1,ymin=-1,xmax=1,ymax=1] {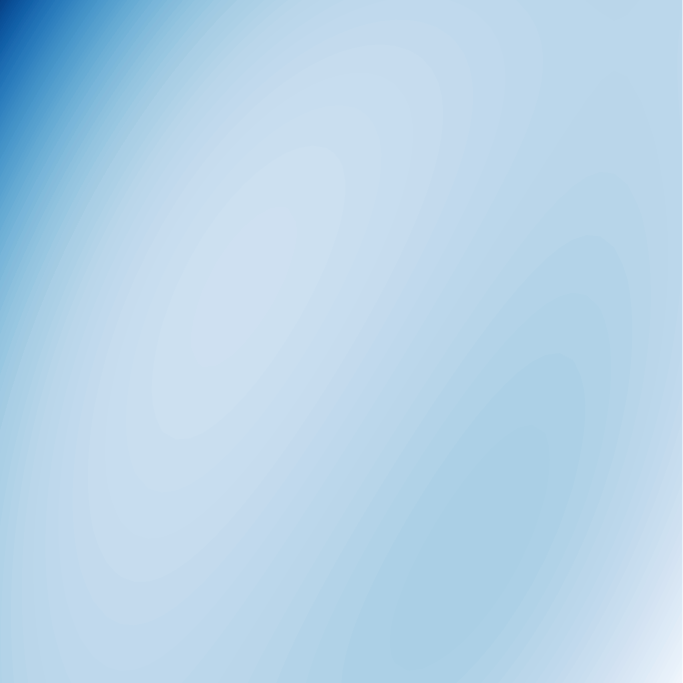};
            
            \draw[black] (axis cs:-1,-1) -- (axis cs:-1,1) -- (axis cs:1,1) -- (axis cs:1,-1) -- cycle;
         
         \draw[fill=Set1-A,Set1-A] (-1,0.86114) circle (0.018);
\draw[fill=black,black] (-1,0.33998) circle (0.018);
\draw[fill=black,black] (-1,-0.33998) circle (0.018);
\draw[fill=black,black] (-1,-0.86114) circle (0.018);
\draw[fill=black,black] (0.86114,1) circle (0.018);
\draw[fill=black,black] (0.33998,1) circle (0.018);
\draw[fill=black,black] (-0.33998,1) circle (0.018);
\draw[fill=black,black] (-0.86114,1) circle (0.018);
\draw[fill=black,black] (1,-0.86114) circle (0.018);
\draw[fill=black,black] (1,-0.33998) circle (0.018);
\draw[fill=black,black] (1,0.33998) circle (0.018);
\draw[fill=black,black] (1,0.86114) circle (0.018);
\draw[fill=black,black] (-0.86114,-1) circle (0.018);
\draw[fill=black,black] (-0.33998,-1) circle (0.018);
\draw[fill=black,black] (0.33998,-1) circle (0.018);
\draw[fill=black,black] (0.86114,-1) circle (0.018);   
         
		\end{axis}
		\begin{scope}[on behind layer]
            \draw[black!00] (0,-0.25) rectangle (5,5);
        \end{scope}
\end{tikzpicture}

%% file: figs/tikz/quad_p3_corr_f2_dg_to.tex
    \begin{tikzpicture}[scale=2]
		\begin{axis}[name=plot1,xlabel={$x$},ylabel={$y$},
    		axis line style={latex-latex},
            axis y line=middle,
            axis x line=middle,
            xmode=linear, % not log
            ymode=linear, % not log
		    xtick={0},
		    ytick={0},
		    xticklabels={ },
		    yticklabels={ },
		    ylabel style={rotate=0, shift={(-0.09,0.06)}},
		    xlabel style={rotate=0, shift={(0.03,0.01)}},
    		xmin=-1.1,xmax=1.1,
    		ymin=-1.1,ymax=1.1,
    		width=5cm,
    		height=5cm,
			scale only axis=true,
    		style={font=\normalsize},
    		axis on top,
    		colorbar,
    		colormap/Blues,
    		colorbar style={
    		    point meta max=5, point meta min=-2,
                ylabel=$\nabla\cdot\mathbf{h}$,
                ylabel style={font=\normalsize},
                ytick={-2,0,2,4},
                yticklabel style={
                    text width=2em,
                    align=left,
                    font=\normalsize,
                }
            }
        ]
            \addplot[thick] graphics[xmin=-1,ymin=-1,xmax=1,ymax=1] {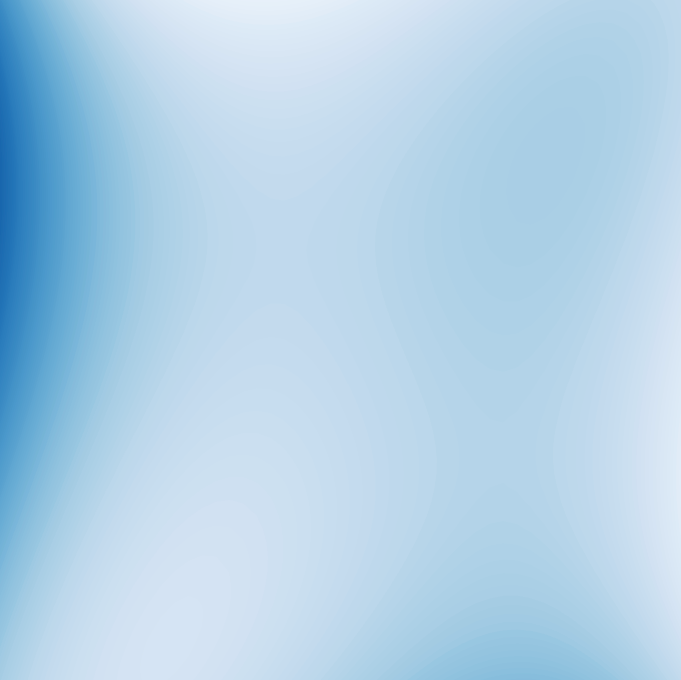};
            
            \draw[black] (axis cs:-1,-1) -- (axis cs:-1,1) -- (axis cs:1,1) -- (axis cs:1,-1) -- cycle;
            
            \draw[fill=black,black] (-1,0.86114) circle (0.018);
\draw[fill=Set1-A,Set1-A] (-1,0.33998) circle (0.018);
\draw[fill=black,black] (-1,-0.33998) circle (0.018);
\draw[fill=black,black] (-1,-0.86114) circle (0.018);
\draw[fill=black,black] (0.86114,1) circle (0.018);
\draw[fill=black,black] (0.33998,1) circle (0.018);
\draw[fill=black,black] (-0.33998,1) circle (0.018);
\draw[fill=black,black] (-0.86114,1) circle (0.018);
\draw[fill=black,black] (1,-0.86114) circle (0.018);
\draw[fill=black,black] (1,-0.33998) circle (0.018);
\draw[fill=black,black] (1,0.33998) circle (0.018);
\draw[fill=black,black] (1,0.86114) circle (0.018);
\draw[fill=black,black] (-0.86114,-1) circle (0.018);
\draw[fill=black,black] (-0.33998,-1) circle (0.018);
\draw[fill=black,black] (0.33998,-1) circle (0.018);
\draw[fill=black,black] (0.86114,-1) circle (0.018);  
         
		\end{axis}
\end{tikzpicture}

%% file: erfr_quad_eo_near.tex
\section{Approximate Euclidean order basis}\label{sec:quad_neo}
    Across two works \citep{Trefethen2017,Trefethen2017b}, Trefethen investigated a Euclidean basis where $\|\mathbf{k}\|_2 \leq k_\mathrm{max}$. In these works a paradox is pointed out: the total order basis is isotropic in the sense that the orders in various directions are equal, however the hyper-cube is exponentially anisotropic, and functions typically require higher orders along diagonals. The conclusion is that a truly isotropic basis for a hyper-cube is more similar to a Euclidean basis. Clearly, from the perspective of applications such as resolving turbulent features within a flow field, we would like the numerical properties to be as isotropic as possible. Therefore, here we consider defining the flux reconstruction scheme on a Euclidean basis.
    
    As discussed in \cref{sec:quad_to}, the symmetry orbits of a quadrilateral place a limit on the set of solution points. For a total order basis we can avoid this problem with the Padua points, as they are provably optimal in some respects; however, no analogous point set currently exists for a Euclidean basis. Therefore, a reasonable alternative is to increase the number of basis functions slightly so that they correspond to a number of points that can be found within the orbits of a quadrilateral. To do this, we can increase $p$ in $\|\mathbf{k}\|_p\leq k_\mathrm{max}$ until a symmetrical set of orbits can be found. This does not need to be performed with any great accuracy due to the discrete nature of the problem. 
    
    \begin{table}[tbhp]
        \centering
        \begin{tabular}{l | r r r r r r r r r r }
            \toprule
            $k_\mathrm{max}$ & 2 & 3 & 4 & 5 & 6 & 7 & 8 & 9 & 10 & 11\\ \midrule
            $p$ & 48 & 50 & 2 & 21 & 3 & 3 & 2.2 & 2 & 2.4 & 2.2 \\
            $n_b$ & 8 & 13 & 17 & 29 & 37 & 45 & 60 & 73 & 92 & 109\\\bottomrule
        \end{tabular}
        \caption{\label{tab:near_eo}Approximate Euclidean basis $p$ and $n_b$ for various orders.}
    \end{table}
    
    \cref{tab:near_eo} shows the approximate values of $p$ and $n_b$ for various orders. We will call this basis an \emph{approximate} Euclidean basis and we use the notation of $p=2^*$ to indicate this. We now enumerate the resulting FR $\tb{Q}$ matrices and stability conditions for several of these orders.
    
\subsection{$k=2$}
    Unlike the true Euclidean order basis at $k=2$, the approximate Euclidean basis has more points than the total order basis. We find that:
    \begin{equation}
        \tb{Q} = \begin{bmatrix}
            0 & 0 & 0 & 0 & 0 & 0 & 0 & 0 \\
            0 & 0 & 0 & 0 & 0 & 0 & 0 & -3q_1 \\
            0 & 0 & 0 & 0 & 0 & 0 & 9q_1 & 0 \\
            0 & 0 & 0 & 0 & 0 & -3q_1 & 0 & 0 \\
            0 & 0 & 0 & 0 & q_1 & 0 & 0 & 0 \\
            0 & 0 & 0 & -3q_1 & 0 & q_0 & 0 & 0 \\
            0 & 0 & 9q_1 & 0 & 0 & 0 & 0 & 0 \\
            0 & -3q_1 & 0 & 0 & 0 & 0 & 0 & q_0
        \end{bmatrix}, \quad \mathrm{for}\quad \mathbf{\Phi}_2 = \begin{bmatrix}
            1\\
            \psi_0(x)\psi_1(y)\\
            \psi_0(x)\psi_2(y)\\
            \psi_1(x)\psi_0(y)\\
            \psi_1(x)\psi_1(y)\\
            \psi_1(x)\psi_2(y)\\
            \psi_2(x)\psi_0(y)\\
            \psi_2(x)\psi_1(y)
        \end{bmatrix}.
    \end{equation}
    This is subject to the stability conditions stemming from positive definiteness and leads to the inequalities:
    \begin{equation}
        60q_0 - 405q_1^2 > -16 \quad \mathrm{and} \quad 2025q_1^2 < 16.
    \end{equation}

\subsection{$k=3$}
    Repeating this for $k=3$ we find: 
    \begin{equation}
        \tb{Q} = \left[ \begin{array}{ccccccc}
            \mathbf{0} & \multicolumn{6}{c}{\mathbf{0}} \\
            \multirow{6}{*}{$\boldsymbol{0}$} & q_0 & 0 & 0 & 0 & 0 & q_2 \\
             & 0 & 0 & 0 & 0 & 0 & 0 \\
             & 0 & 0 & 0 & 0 & 0 & 0 \\
             & 0 & 0 & 0 & q_1 & 0 & 0 \\
             & 0 & 0 & 0 & 0 & 0 & 0 \\
             & q_2 & 0 & 0 & 0 & 0 & q_0 \\
        \end{array}\right], \quad \mathrm{for}\quad \mathbf{\Phi}_3=\begin{bmatrix}
            1\\
            \vdots\\
            \psi_1(x)\psi_3(y)\\
            \psi_2(x)\psi_0(y)\\
            \psi_2(x)\psi_1(y)\\
            \psi_2(x)\psi_2(y)\\
            \psi_3(x)\psi_0(y)\\
            \psi_3(x)\psi_1(y)\\
        \end{bmatrix}.
    \end{equation}
    which is subject to the stability conditions that:
    \begin{equation}
        21q_0 > -4, \quad 25q_1 > -4, \quad \mathrm{and} \quad21q_0(21q_0 + 8) - 441q_2^2 > -16.
    \end{equation}
    
    An example of an approximate Euclidean order basis correction function is included in \cref{fig:quad_neo_k3_dg} for $\tb{Q}=0$. Comparison with the correction function shown in \cref{fig:quad_to_k3_dg} shows subtle differences, most notably in the ranges of the respective functions.
    \begin{figure}[tbhp]
        \centering
        \subfloat[First flux point $(-1,\frac{\sqrt{15+2\sqrt{30}}}{35})$.]{\adjustbox{width=0.39\linewidth,valign=b}{\input{figs/tikz/quad_k3_corr_f1_dg_neo}}}
        ~
        \subfloat[Second flux point $(-1,\frac{\sqrt{15-2\sqrt{30}}}{35})$.]{\adjustbox{width=0.56\linewidth,valign=b}{\input{figs/tikz/quad_k3_corr_f2_dg_neo}}}
        \caption{\label{fig:quad_neo_k3_dg}Divergence of DG correction field for $k=3$ FR on a quadrilateral with an approximate Euclidean order basis for two flux points, shown in red.}
    \end{figure}

\subsection{$k=4$}\label{ssec:quad_eo_k4}
    We can repeat this analysis again for $k=4$; however, in this case the approximate Euclidean order basis and the Euclidean order basis are the same. We then find that:
    \begin{equation}
        \tb{Q} = \left[ \begin{array}{cccccccccccccc}
            \mathbf{0}  & \multicolumn{13}{c}{\mathbf{0}} \\
            \multirow{13}{*}{$\boldsymbol{0}$}  & q_0 & 0 & 0 & 0 & 0 & 0 & 0 & 0 & 0 & 0 & 0 & 0 & q_3 \\
          & 0 & 0 & 0 & 0 & 0 & 0 & 0 & 0 & 0 & 0 & 0 & 0 & 0 \\
          & 0 & 0 & 0 & 0 & 0 & 0 & 0 & 0 & 0 & 0 & 0 & 0 & 0 \\
          & 0 & 0 & 0 & 0 & 0 & 0 & 0 & 0 & 0 & 0 & 0 & -5q_2/3 & 0 \\
          & 0 & 0 & 0 & 0 & 0 & 0 & 0 & 0 & 0 & 0 & 25q_2/9 & 0 & 0 \\
          & 0 & 0 & 0 & 0 & 0 & 0 & 0 & 0 & 0 & 0 & 0 & 0 & 0 \\
          & 0 & 0 & 0 & 0 & 0 & 0 & 0 & 0 & -5q_2/3 & 0 & 0 & 0 & 0 \\
          & 0 & 0 & 0 & 0 & 0 & 0 & 0 & q_2 & 0 & 0 & 0 & 0 & 0 \\
          & 0 & 0 & 0 & 0 & 0 & 0 & -5q_2/3 & 0 & q_1 & 0 & 0 & 0 & 0 \\
          & 0 & 0 & 0 & 0 & 0 & 0 & 0 & 0 & 0 & 0 & 0 & 0 & 0 \\
          & 0 & 0 & 0 & 0 & 25q_2/9 & 0 & 0 & 0 & 0 & 0 & 0 & 0 & 0 \\
          & 0 & 0 & 0 & -5q_2/3 & 0 & 0 & 0 & 0 & 0 & 0 & 0 & q_1 & 0 \\
          & q_3 & 0 & 0 & 0 & 0 & 0 & 0 & 0 & 0 & 0 & 0 & 0 & q_0
        \end{array}\right]%, \quad \mathrm{for}\quad \mathbf{\Phi}_4=\begin{bmatrix}
        %     1 \\
        %     \vdots \\
        %     \psi_0(x)\psi_4(y) \\
        %     \psi_1(x)\psi_0(y) \\
        %     \psi_1(x)\psi_1(y) \\
        %     \psi_1(x)\psi_2(y) \\
        %     \psi_1(x)\psi_3(y) \\
        %     \psi_2(x)\psi_0(y) \\
        %     \psi_2(x)\psi_1(y) \\
        %     \psi_2(x)\psi_2(y) \\
        %     \psi_2(x)\psi_3(y) \\
        %     \psi_3(x)\psi_0(y) \\
        %     \psi_3(x)\psi_1(y) \\
        %     \psi_3(x)\psi_2(y) \\
        %     \psi_4(x)\psi_0(y)
        % \end{bmatrix}
    \end{equation}
    subject to the stability constraints that:
    \begin{subequations}
        \begin{align}
            q_0 &>-\frac{4}{9},\\
            q_2 &> -\frac{4}{25},\\
            420q_1 + 48 - 4375q_2^2 &> 0,\\
            144 - 30625q_2^2 &> 0,\\
            9q_0(9q_0 + 8) - 81q_3^2 + 16 &> 0.
        \end{align}
    \end{subequations}

%% file: figs/tikz/quad_k3_corr_f1_dg_neo.tex
    \begin{tikzpicture}[scale=2]
		\begin{axis}[name=plot1,xlabel={$x$},ylabel={$y$},
    		axis line style={latex-latex},
            axis y line=middle,
            axis x line=middle,
            xmode=linear, % not log
            ymode=linear, % not log
		    xtick={0},
		    ytick={0},
		    xticklabels={ },
		    yticklabels={ },
		    ylabel style={rotate=0, shift={(-0.09,0.06)}},
		    xlabel style={rotate=0, shift={(0.03,0.01)}},
    		xmin=-1.1,xmax=1.1,
    		ymin=-1.1,ymax=1.1,
    		width=5cm,
    		height=5cm,
			scale only axis=true,
    		style={font=\normalsize},
    		axis on top,
    % 		colorbar,
    % 		colormap/Blues,
    % 		colorbar style={
    % 		    point meta max=5, point meta min=-2,
    %             ylabel=$\nabla\cdot\mathbf{h}$,
    %             ylabel style={font=\normalsize},
    %             ytick={-4,-2,0,2},
    %             yticklabel style={
    %                 text width=2em,
    %                 align=left,
    %                 font=\normalsize,
    %             }
    %         }
        ]
            \addplot[thick] graphics[xmin=-1,ymin=-1,xmax=1,ymax=1] {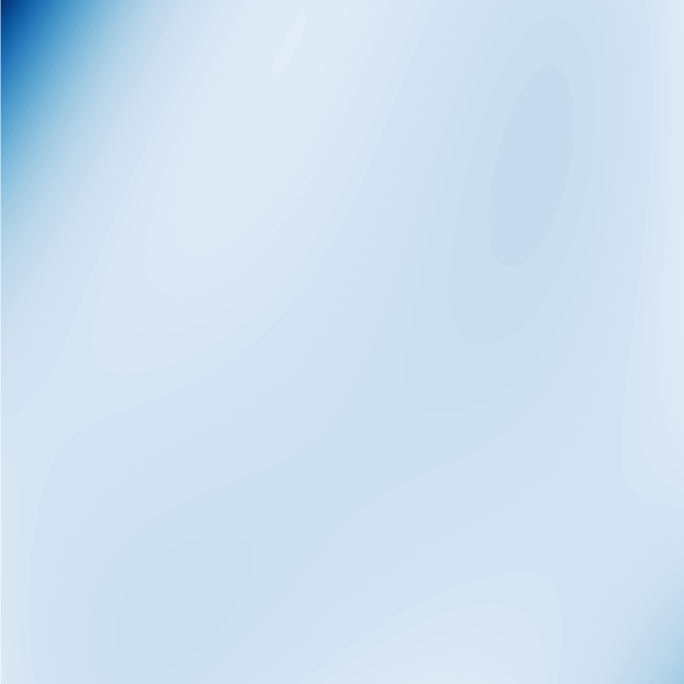};
            
            \draw[black] (axis cs:-1,-1) -- (axis cs:-1,1) -- (axis cs:1,1) -- (axis cs:1,-1) -- cycle;
         
         \draw[fill=Set1-A,Set1-A] (-1,0.86114) circle (0.018);
\draw[fill=black,black] (-1,0.33998) circle (0.018);
\draw[fill=black,black] (-1,-0.33998) circle (0.018);
\draw[fill=black,black] (-1,-0.86114) circle (0.018);
\draw[fill=black,black] (0.86114,1) circle (0.018);
\draw[fill=black,black] (0.33998,1) circle (0.018);
\draw[fill=black,black] (-0.33998,1) circle (0.018);
\draw[fill=black,black] (-0.86114,1) circle (0.018);
\draw[fill=black,black] (1,-0.86114) circle (0.018);
\draw[fill=black,black] (1,-0.33998) circle (0.018);
\draw[fill=black,black] (1,0.33998) circle (0.018);
\draw[fill=black,black] (1,0.86114) circle (0.018);
\draw[fill=black,black] (-0.86114,-1) circle (0.018);
\draw[fill=black,black] (-0.33998,-1) circle (0.018);
\draw[fill=black,black] (0.33998,-1) circle (0.018);
\draw[fill=black,black] (0.86114,-1) circle (0.018);   
         
		\end{axis}
		\begin{scope}[on behind layer]
            \draw[black!00] (0,-0.25) rectangle (5,5);
        \end{scope}
\end{tikzpicture}

%% file: figs/tikz/quad_k3_corr_f2_dg_neo.tex
    \begin{tikzpicture}[scale=2]
		\begin{axis}[name=plot1,xlabel={$x$},ylabel={$y$},
    		axis line style={latex-latex},
            axis y line=middle,
            axis x line=middle,
            xmode=linear, % not log
            ymode=linear, % not log
		    xtick={0},
		    ytick={0},
		    xticklabels={ },
		    yticklabels={ },
		    ylabel style={rotate=0, shift={(-0.09,0.06)}},
		    xlabel style={rotate=0, shift={(0.03,0.01)}},
    		xmin=-1.1,xmax=1.1,
    		ymin=-1.1,ymax=1.1,
    		width=5cm,
    		height=5cm,
			scale only axis=true,
    		style={font=\normalsize},
    		axis on top,
    		colorbar,
    		colormap/Blues,
    		colorbar style={
    		    point meta max=8, point meta min=-2,
                ylabel=$\nabla\cdot\mathbf{h}$,
                ylabel style={font=\normalsize},
                ytick={-2,0,2,4,6,8},
                yticklabel style={
                    text width=2em,
                    align=left,
                    font=\normalsize,
                }
            }
        ]
            \addplot[thick] graphics[xmin=-1,ymin=-1,xmax=1,ymax=1] {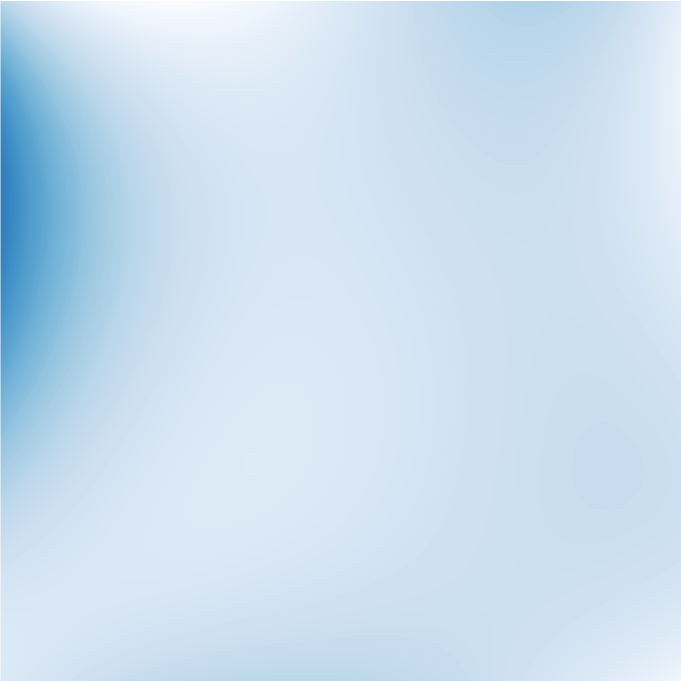};
            
            \draw[black] (axis cs:-1,-1) -- (axis cs:-1,1) -- (axis cs:1,1) -- (axis cs:1,-1) -- cycle;
            
            \draw[fill=black,black] (-1,0.86114) circle (0.018);
            
            \draw[fill=Set1-A,Set1-A] (-1,0.33998) circle (0.018);
            \draw[fill=black,black] (-1,-0.33998) circle (0.018);
            \draw[fill=black,black] (-1,-0.86114) circle (0.018);
            \draw[fill=black,black] (0.86114,1) circle (0.018);
            \draw[fill=black,black] (0.33998,1) circle (0.018);
            \draw[fill=black,black] (-0.33998,1) circle (0.018);
            \draw[fill=black,black] (-0.86114,1) circle (0.018);
            \draw[fill=black,black] (1,-0.86114) circle (0.018);
            \draw[fill=black,black] (1,-0.33998) circle (0.018);
            \draw[fill=black,black] (1,0.33998) circle (0.018);
            \draw[fill=black,black] (1,0.86114) circle (0.018);
            \draw[fill=black,black] (-0.86114,-1) circle (0.018);
            \draw[fill=black,black] (-0.33998,-1) circle (0.018);
            \draw[fill=black,black] (0.33998,-1) circle (0.018);
            \draw[fill=black,black] (0.86114,-1) circle (0.018);  
         
		\end{axis}
\end{tikzpicture}

%% file: numerical.tex
\section{Numerical Experiments}\label{sec:numerical}
    In this section we present results of numerical experiments with the linear advection equation. In particular, we are concerned with: 
    \begin{equation}
        \px{u}{t} + \nabla\cdot\mathbf{a}u, \quad \mathrm{for}\quad \mathbf{a}=\begin{bmatrix}
            \cos\theta \\\sin\theta
        \end{bmatrix}.
    \end{equation}
    To test the effects of anisotropy, we use an initial condition comprised of several superimposed Morlet wavelets~\citep{Kronland1987}, with the definition:
    \begin{subequations}
        \begin{align}
            u &= c_\sigma\pi^{-1/4}\sum^n_{i=1}\exp(-r_i^2/2)\left[\cos(\sigma r_i) - \kappa_i\right],\\
            c_\sigma &= \left[1 + \exp(-\sigma^2) - 2\exp(-3\sigma^2/4)\right]^{-1/2},\\
            r_i &= \sqrt{(x-x_i)^2 + (y-y_i)^2},
        \end{align}
    \end{subequations}
    where $(x_i,y_i)$ is a random centre coordinate, and $\sigma$ and $\kappa_i$ are control parameters. For the experiments conducted, four wavelets were superimposed, $n=4$, with the control parameter $\sigma$ set to three and $\kappa_i\in[0,1]$ randomly chosen for each wavelet. A series of advection angles were tested and the initial condition for each was the same, with the same random numbers chosen via a Mersenne twister algorithm. This initial condition is ideal for testing isotropy due to the dependence on radius and wider frequency spectrum.
    
    The domain used was fully periodic and covered $K\in[0,2\pi]^2$, partitioned into $N$ regular quadrilaterals. For time integration an explicit SSP-RK3 scheme was used with constant $\Delta t=10^{-3}$, and for all tests the common interfaces were fully upwinded.

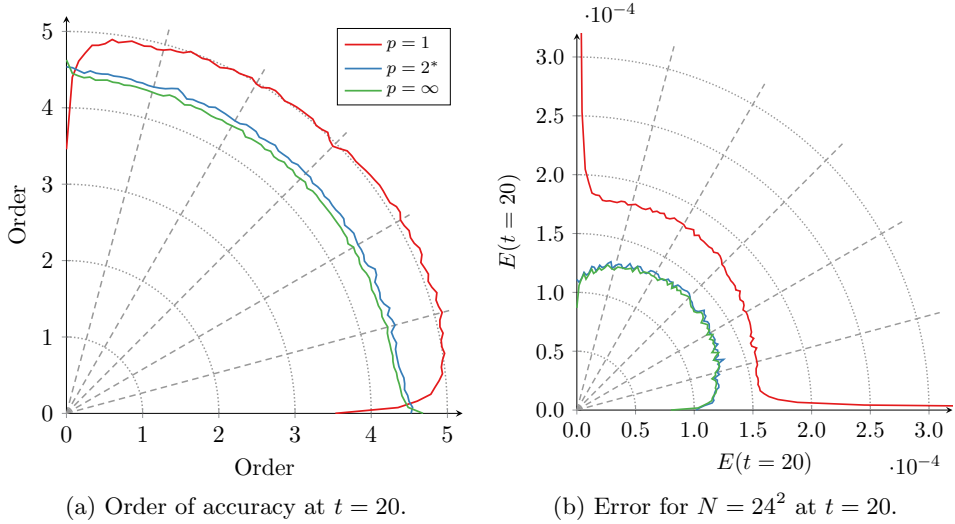
\begin{figure}[tbhp]
        \centering
        \subfloat[Order of accuracy at $t=20$.]{\label{fig:k3_wm_order}\adjustbox{width=0.4\linewidth,valign=b}{\input{figs/tikz/order_circ_mw_k3_waves}}}
        ~
        \subfloat[Error for $N=24^2$ at $t=20$.]{\label{fig:k3_mw_error}\adjustbox{width=0.4\linewidth,valign=b}{\input{figs/tikz/error_circ_mw_k3_waves}}}
        \caption{\label{fig:k3_iso_mw}Variation of order and error with angle, $\theta$, for $k=3$ DG FR with different bases.}
    \end{figure}

    Initially, a sweep of advection angles for $N\in\{8^2,10^2,\dots,32^2\}$ was performed, the results of which are presented in \cref{fig:k3_iso_mw}. This shows a marked difference between the total order, approximate Euclidean order ($p=2^*$), and maximal order bases. Most notably, the error when using a total order basis is significantly higher. This is consistent with the findings of \citet{Trefethen2017} for the interpolation error of the two-dimensional Runge function.
    
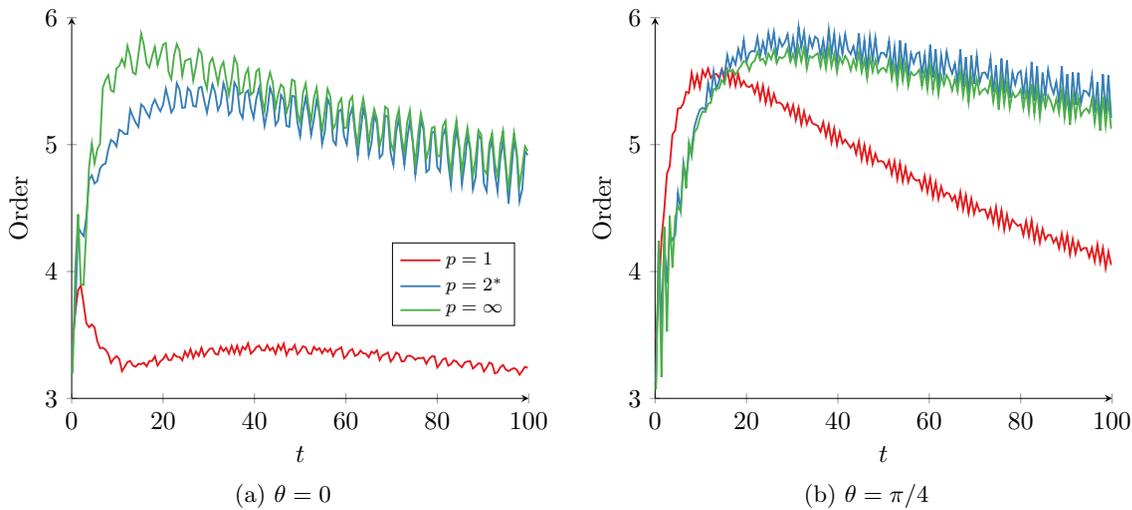
\begin{figure}[tbhp]
        \centering
        \subfloat[$\theta=0$]{\label{fig:mw_k3_order_t_0}\adjustbox{width=0.48\linewidth,valign=b}{\input{figs/tikz/order_mw_k3_theta0}}}
        ~
        \subfloat[$\theta=\pi/4$]{\adjustbox{width=0.48\linewidth,valign=b}{\input{figs/tikz/order_mw_k3_theta45}}}
        \caption{\label{fig:morelt_order}Order versus time for $k=3$ DG FR with different bases, calculated for $N=\{8^2,12^2\}$.}
    \end{figure}

    From \cref{fig:k3_wm_order}, we see that the order of accuracy of the total order scheme is higher for a large range of angles. Investigating this further, we present the variation of order in time calculated for two grids ($N\in\{8^2,12^2\}$) and two angles, see \cref{fig:morelt_order}. This shows that for non-grid aligned angles, the decay of the low order secondary modes is faster, seen by the faster transition from order $k+1$ to $2k$. This is responsible for the apparently higher order shown in \cref{fig:k3_wm_order}. However, after the peak order of $2k$ is reached~\citep{Asthana2017}, the decay towards order $k+1$, is faster and is generally indicative of the total order basis having larger dispersion and diffusion errors at higher frequencies. Decay in the order is seen for all bases as time progresses, and is due to dispersion errors at high frequencies. A further effect of the total order basis is observed in  \cref{fig:mw_k3_order_t_0}, where for grid aligned waves the total order basis does not exhibit the super-convergence property observed for the other bases. Furthermore, from \cref{fig:k3_mw_error} it is clear that the error when using a total order basis is asymmetric about $\theta=\pi/4$, with lower error observed at $\theta=\pi/6$ than $\theta=\pi/3$. Given that this is not found to occur for the other bases, this is a direct result of the Padua points lacking full rotational symmetry on a quadrilateral, as can be seen in \cref{fig:padua}.
    
    Finally, the points used for the $p=2^*$ cases were optimised to reduce $L_2$ error, as this has previously been shown to be important in one dimension~\citep{Witherden2015}. As alternatives, Lebesgue and Fekete optimal point sets were also produced --- the results of which are not shown here, but which were significantly worse than those with $L_2$ optimised points in terms of absolute error.

%% file: figs/tikz/order_circ_mw_k3_waves.tex
\begin{tikzpicture}
    \definecolor{agrey}{rgb}{0.6,0.6,0.6}

    \begin{axis}
    [
        width=8cm,
        height=8cm,
        axis line style={latex-latex},
        axis y line=left,
        axis x line=left,
        xmode=linear, % not log
        ymode=linear, % not log
        xlabel = {Order},
        ylabel = {Order},
        xmin = 0, xmax = 5.2,
        ymin = 0, ymax = 5.2,
        legend cell align={left},
        legend style={font=\scriptsize, at={(0.97, 0.97)},anchor=north east, 
        %/tikz/column 2/.style={column sep=5pt}
        },
        %axis line style={draw=none},
        %tick style={draw=none},
        x tick label style={/pgf/number format/.cd, fixed, fixed zerofill, precision=0, /tikz/.cd},
        y tick label style={/pgf/number format/.cd, fixed, fixed zerofill, precision=0, /tikz/.cd},
    ]
        \foreach \t in {15,30,45,60,75} 
        {
            \addplot[forget plot, color=agrey, line width=0.7pt, densely dashed] coordinates{(0, 0) ({5.2*cos(\t)}, {5.2*sin(\t)})};
        }
        
        \foreach \r in {1, 2, 3, 4, 5} 
        {
            \addplot[forget plot, color=agrey, samples=40, domain=0.0*pi:0.5*pi, line width=0.7pt, densely dotted]({\r*cos(deg(x))}, {\r*sin(deg(x))});
        }
    
        \addplot[thick, color={Set1-A}] table[x expr={(\thisrowno{1})*cos(\thisrowno{0})}, y expr={(\thisrowno{1})*sin(\thisrowno{0})}, col sep=comma]{./figs/data/morlet_odata_k3_lp1.csv};
                                                  
        \addplot[thick, color={Set1-B}] table[x expr={(\thisrowno{1})*cos(\thisrowno{0})}, y expr={(\thisrowno{1})*sin(\thisrowno{0})}, col sep=comma]{./figs/data/morlet_odata_rob_k3_lp50.csv};
        
        \addplot[thick, color={Set1-C}] table[x expr={(\thisrowno{1})*cos(\thisrowno{0})}, y expr={(\thisrowno{1})*sin(\thisrowno{0})}, col sep=comma]{./figs/data/morlet_odata_k3_lpInf.csv};

        \addlegendentry{$p=1$};
        \addlegendentry{$p=2^*$};
        \addlegendentry{$p=\infty$};
  
    \end{axis}
    
\end{tikzpicture}

%% file: figs/tikz/error_circ_mw_k3_waves.tex
\begin{tikzpicture}
    \definecolor{agrey}{rgb}{0.6,0.6,0.6}

    \begin{axis}
    [
        width=8cm,
        height=8cm,
        axis line style={latex-latex},
        axis y line=left,
        axis x line=left,
        xmode=linear, % not log
        ymode=linear, % not log
        xlabel = {$E(t=20)$},
        ylabel = {$E(t=20)$},
        xmin = 0, xmax = 3.2e-4,
        ymin = 0, ymax = 3.2e-4,
        legend cell align={left},
        legend style={font=\scriptsize, at={(0.97, 0.97)},anchor=north east, 
        %/tikz/column 2/.style={column sep=5pt}
        },
        %axis line style={draw=none},
        %tick style={draw=none},
        x tick label style={/pgf/number format/.cd, fixed, fixed zerofill, precision=1, /tikz/.cd},
        y tick label style={/pgf/number format/.cd, fixed, fixed zerofill, precision=1, /tikz/.cd},
    ]
        \foreach \t in {15,30,45,60,75} 
        {
            \addplot[forget plot, color=agrey, line width=0.7pt, densely dashed] coordinates{(0, 0) ({3.2e-4*cos(\t)}, {3.2e-4*sin(\t)})};
        }
        
        \foreach \r in {0.5e-4,1e-4, 1.5e-4,2e-4, 2.5e-4, 3e-4} 
        {
            \addplot[forget plot, color=agrey, samples=40, domain=0.0*pi:0.5*pi, line width=0.7pt, densely dotted]({\r*cos(deg(x))}, {\r*sin(deg(x))});
        }
    
        \addplot[thick, color={Set1-A}] table[x expr={(\thisrowno{2})*cos(\thisrowno{0})}, y expr={(\thisrowno{2})*sin(\thisrowno{0})}, col sep=comma]{./figs/data/morlet_odata_k3_lp1.csv};
                                                  
        \addplot[thick, color={Set1-B}] table[x expr={(\thisrowno{2})*cos(\thisrowno{0})}, y expr={(\thisrowno{2})*sin(\thisrowno{0})}, col sep=comma]{./figs/data/morlet_odata_rob_k3_lp50.csv};
        
        \addplot[thick, color={Set1-C}] table[x expr={(\thisrowno{2})*cos(\thisrowno{0})}, y expr={(\thisrowno{2})*sin(\thisrowno{0})}, col sep=comma]{./figs/data/morlet_odata_k3_lpInf.csv};

        % \addlegendentry{$p=1$};
        % \addlegendentry{$p=2^*$}
        % \addlegendentry{$p=\infty$};
  
    \end{axis}
    
\end{tikzpicture}

%% file: figs/tikz/order_mw_k3_theta0.tex
\begin{tikzpicture}
    \begin{axis}
    [
        axis line style={latex-latex},
        axis y line=left,
        axis x line=left,
        xmode=linear, % not log
        ymode=linear, % not log
        xlabel = {$t$},
        ylabel = {Order},
        xmin = 0, xmax = 100,
        ymin = 3, ymax = 6,
        ytick = {3,4,5,6},
        legend cell align={left},
        legend style={font=\scriptsize, at={(0.97, 0.3)},anchor=east},
        %axis line style={draw=none},
        %tick style={draw=none},
        x tick label style={/pgf/number format/.cd, fixed, fixed zerofill, precision=0, /tikz/.cd},
        y tick label style={/pgf/number format/.cd, fixed, fixed zerofill, precision=0, /tikz/.cd},
    ]
        
        \addplot[color=Set1-A,thick] table[x=t, y=o, col sep=comma]{./figs/data/order_mw_k3_lp1_theta0.csv};
        
        \addplot[color=Set1-B,thick] table[x=t, y=o, col sep=comma]{./figs/data/order_mw_k3_lp50_theta0.csv};
        
        \addplot[color=Set1-C,thick] table[x=t, y=o, col sep=comma]{./figs/data/order_mw_k3_lpInf_theta0.csv};
        
        \addlegendentry{$p=1$};
        \addlegendentry{$p=2^*$};
        \addlegendentry{$p=\infty$};
        
    \end{axis}
\end{tikzpicture}

%% file: figs/tikz/order_mw_k3_theta45.tex
\begin{tikzpicture}
    \begin{axis}
    [
        axis line style={latex-latex},
        axis y line=left,
        axis x line=left,
        xmode=linear, % not log
        ymode=linear, % not log
        xlabel = {$t$},
        ylabel = {Order},
        xmin = 0, xmax = 100,
        ymin = 3, ymax = 6,
        ytick = {3,4,5,6},
        legend cell align={left},
        legend style={font=\scriptsize, at={(0.97, 0.3)},anchor=east},
        %axis line style={draw=none},
        %tick style={draw=none},
        x tick label style={/pgf/number format/.cd, fixed, fixed zerofill, precision=0, /tikz/.cd},
        y tick label style={/pgf/number format/.cd, fixed, fixed zerofill, precision=0, /tikz/.cd},
    ]
        
        \addplot[color=Set1-A,thick] table[x=t, y=o, col sep=comma]{./figs/data/order_mw_k3_lp1_theta45.csv};
        
        \addplot[color=Set1-B,thick] table[x=t, y=o, col sep=comma]{./figs/data/order_mw_k3_lp50_theta45.csv};
        
        \addplot[color=Set1-C,thick] table[x=t, y=o, col sep=comma]{./figs/data/order_mw_k3_lpInf_theta45.csv};
        
        %\addlegendentry{$p=1$};
        %\addlegendentry{$p=2^*$};
        %\addlegendentry{$p=\infty$};
        
    \end{axis}
\end{tikzpicture}